\documentclass[10pt]{elsarticle}

\usepackage{amsmath,amssymb,amsfonts,color}
\usepackage{algorithmic}
\usepackage{algorithm,booktabs}
\usepackage[mathscr]{eucal}

\DeclareMathOperator*{\argmin}{arg\rm{}min}

\renewcommand{\theequation}{\arabic{section}.\arabic{equation}}

\newcommand{\R}{\mathbb{R}}

\journal{Advances in Computational Mathematics}

\begin{document}

\begin{frontmatter}



\title{Randomized Model Order Reduction} 
%via\\ Proper Orthogonal Decomposition and \\Dynamic Mode Decomposition}


\author{Alessandro Alla\corref{Thanks1}}
\ead{aalla@fsu.edu}
\cortext[Thanks1]{This author wishes to acknowledge the support of the Department of Energy grant number DE-SC0009324.}
\address{Department of Scientific Computing,,Florida State University, Tallahassee FL-32306}
\author{J. Nathan Kutz \corref{Thanks2}}
\ead{kutz@uw.edu}
\cortext[Thanks2]{This author gratefully acknowledges support from the U.S. Air Force Office of Scientific Research (FA9550-15-1-0385).}
\address{Department of Applied Mathematics, University of Washington, Seattle, WA 98195-3925}

\begin{abstract}
Singular value decomposition (SVD) has a crucial role in model order reduction. It is often utilized in the offline stage to compute basis functions that project the high-dimensional nonlinear problem into a low-dimensionsl model which is, then, evaluated cheaply. It constitutes a building block for many techniques such as e.g. Proper Orthogonal Decomposition and Dynamic Mode Decomposition.\\
The aim of this work is to provide efficient computation of the basis functions via randomized matrix decompositions. This is possible due to the randomized Singular Value Decomposition (rSVD) which is a fast and accurate alternative of the SVD. Although this is considered as offline stage, this computation may be extremely expensive and therefore the use of compressed techniques drastically reduce its cost.  Numerical examples show the effectiveness of the method for both Proper Orthogonal Decomposition (POD) and Dynamic Mode Decomposition (DMD).
 \end{abstract}


\end{frontmatter}

%\maketitle

%%%%%%%%%%%%%%%%%%%%%%%%%%%%%%%%%%%%%%%%%%%%%%%%%%%%%
\section{Introduction}
\label{Section1}
\setcounter{section}{1}
\setcounter{equation}{0}
\setcounter{theorem}{0}
\setcounter{algorithm}{0}
\renewcommand{\theequation}{\arabic{section}.\arabic{equation}}

Reduced order models (ROMs) continue to play a critically enabling role in 
modern, large-scale scientific computing applications~\cite{karen1}. The ROM architecture is
being exploited in many simulation based physics and engineering systems in order
to render tractable many high-dimensional simulations.  Fundamentally,
the ROM algorithmic structure is designed to construct low-dimensional subspaces, typically
computed with singular value decomposition (SVD), where the evolution dynamics can be
embedded using Galerkin projection.  Thus instead of solving a high-dimensional system of
differential equations (e.g. millions or billions of degrees of freedom), a rank $\ell$ model can be
constructed in a principled way.  Three steps are required for this low-rank approximation:  (i) numerical
solutions of the original high-dimensional system, (ii) dimensionality-reduction of this solution data typically 
produced with an SVD, and (iii) Galerkin projection of the dynamics on the low-rank subspace.
The first two steps are often called the {\em offline} stage of the ROM architecture whereas the 
third step is known as the {\em online} stage.  Offline stages are exceptionally expensive, but enable
the (cheap) online stage to potentially run in real time.   In this manuscript, we integrate recent innovations
in randomized linear algebra methods~\cite{gunnar_rLA}, particularly as it relates to the singular value decomposition, and 
compressive sampling in order to (i) improve the computational
efficiency of the second step of the ROM architecture, namely the building of low-rank subspaces
used for Galerkin projection, and (ii) provide a rapid evaluation of the nonlinear terms in the ROM model using compressive sampling of the dynamic mode decomposition.

Randomized methods for matrix computations provide an efficient computation of low-rank structures in data matrices, which
is a foundational aspect of machine learning and big data applications.  
Such algorithms exploit the fact that the target rank of interest, $\ell$, is significantly smaller than the
high-dimensional data under consideration.  In the case of ROMs, there may only be a couple hundred
modes of interest ($\ell\approx 200$) whereas the numerical solution of the original high-dimensional system
may have millions or billions of degrees of freedom.  ROMs allow one to simulate this system with a
differential equation of dimension $\ell$, thus greatly reducing computational time.
Randomized techniques circumvent the challenge of traditional (deterministic) SVD reduction 
which requires significant memory and processing resources for the high-dimensional data generated from full state simulations.  Randomized techniques are robust, reliable and computationally efficient and can be used to construct a smaller (compressed) matrix, which accurately approximates a high-dimensional data matrix~\citep{Mahoney2011,RandNLA}. There exist several strategies for obtaining the compressed matrix, and using random projections is certainly the most robust off-the-shelf approach.  Randomized algorithms have been in particular studied for computing the near-optimal low-rank singular value decomposition by~\cite{frieze2004fast}, \cite{liberty2007randomized}, \cite{woolfe2008fast} and \cite{Martinsson201147}. The seminal contribution~\cite{halko2011rand} extends and surveys this work. 

In addition to randomized techniques, compressive sampling strategies are of growing interest for matrix computations as they also allow for the approximation of decompositions with few measurements.  Much like randomized algorithms, compressive sampling takes advantage of the inherent sparsity of the spatio-temporal dynamics in an appropriate basis.  Thus the target rank $\ell$ determines the number of sample points required.  Compressive sampling can be used with the dynamic mode decomposition (DMD)~\cite{dmdbook} to enact a {\em compressive DMD}~\cite{brunton:cs} approximation for the Galerkin projected dynamics~\cite{AK16}.   The DMD method is an attractive alternative to the standard POD--Galerkin reduction which also uses sparse sampling through the gappy POD and/or DEIM/EIM architecture.  Ultimately, the low-rank structure inherent in ROMs allows the community to exploit sparse measurements to reconstruct an accurate approximation of the high-dimensional system.  In this work, two new innovations are introduced that leverage our current computational capabilities, namely compressive sampling for enhancing the DMD method for ROMs and randomized singular value decompositions for constructing efficient POD basis elements.  

The structure of the paper is as follows.
 In Section \ref{Section2} we review model reduction techniques based upon the Proper Orthogonal Decomposition (POD) method, the Discrete Empirical Interpolation method (DEIM) and the Dynamic Mode Decomposition (DMD) applied to general nonlinear dynamical systems.
 In Section 3 we highlight innovations of randomized techniques for matrix computations, which is the building block for our new approach. Section \ref{Section4} focuses on the application of compressed matrix decompositions in model order reduction. Finally, numerical tests are presented in Section \ref{Section5}. 
Throughout the paper we use the following notation:  all matrices and vectors are in bold letters. The basis functions are denoted by the matrix $\bf {\bf \Psi}$ with different superscripts denoting how we computed the basis, e.g. ${\bf {\bf \Psi}}^{\mbox{\tiny POD}}$ represents the basis functions from the POD method. The rank of the POD basis functions is $\ell$, the rank of the nonlinear term is $k$, whereas $p$ is the number of measurements utilized in the compressed techniques.

%%%%%%%%%%%%%%%%%%%%%%%%%%%%%%%%%%%%%%%%%%%%%%%%%%%%%

%%%%%%%%%%%%%%%%%%%%%%%%%%%%%%%%%%%%%%%%%%%%%%%%%%%%%
\section{Model Order Reduction Techniques}
\label{Section2}
\setcounter{section}{2}
\setcounter{equation}{0}
\setcounter{theorem}{0}
\setcounter{algorithm}{0}
\renewcommand{\theequation}{\arabic{section}.\arabic{equation}}
%%%%%%%%%%%%%%%%%%%%%%%%%%%%%%%%%%%%%%%%%%%%%%%%%%%%%

We consider the general system of high-dimensional, ordinary differential equations:
\begin{equation}
\left\{ \begin{array}{l}\label{ode}
{\bf M}\dot{{\bf y}}(t)={\bf A}{\bf y}(t)+{\bf f}(t,{\bf y}(t)),\;\; t\in(0,T],\\
{\bf y}(0)={\bf y_0},\\
\end{array} \right.
\end{equation}
where ${\bf y_0}\in\R^n$ is a given initial data, ${\bf M, A}\in \R^{n\times n}$ given matrices and ${\bf f}:[0,T]\times\R^n\rightarrow\R^n$ a continuous function in both arguments and locally Lipschitz-type with respect to the second variable. It is well--known that under these assumptions there exists a unique solution for \eqref{ode}.
 This class of problems arises in a wide range of applications, but especially from the numerical approximation of partial differential equations. In such cases, the dimension of the problem $n$ is the number of spatial grid points used from discretization and it typically is very large. The numerical solution of system \eqref{ode} may be very expensive to compute and therefore it is often useful to simplify the complexity of the problem by means of reduced order models. 
The model reduction approach is based on projecting the nonlinear dynamics onto a low dimensional manifold utilizing projectors that contain information from the full, high-dimensional system. 

Let us assume that we have computed some basis functions $ {\bf \Psi}=\{{\boldsymbol{\psi}}_j\}_{j=1}^\ell\in\R^{n\times \ell}$ of rank $\ell$ for  \eqref{ode}.  We can project the dynamics onto the low-rank basis functions using:
\begin{equation}\label{pod_ans}
{\bf y}(t)\approx {\bf \Psi} {\bf y^\ell}(t),
\end{equation}
where ${\bf y}^\ell(t)$ are functions on $\R^\ell$ and defined on the time interval from $[0,T]$. We note that we are working with a Galerkin-type projection where we consider only a few basis functions whose support is non-local, unlike Finite Element basis functions. The reduced solution ${\bf y}^\ell(t)\in V^\ell\subset V$ where $V^\ell=\mbox{span}\{{\boldsymbol{\psi}}_1,\ldots,{\boldsymbol{\psi}}_\ell\}$.

Inserting the projection assumption \eqref{pod_ans} into the full model \eqref{ode}, and making use of the orthogonality of the basis functions, the reduced model takes the following form: 
\begin{equation}\label{pod_sys}
\left\{\begin{array}{l}
{\bf M}^\ell\dot{{\bf y}}^\ell(t)={\bf A}^\ell {\bf y}^\ell(t)+{\bf {\bf \Psi}}^Tf(t,{\bf {\bf \Psi}} {\bf y}^\ell(t)),\quad t\in(0,T],\\
{\bf y}^\ell(0)={\bf y_0^\ell},
\end{array}\right.
\end{equation}
where $({\bf M}^\ell)_{ij}=\langle {\bf M}{\boldsymbol{\psi}}_i,{\boldsymbol{\psi}}_j\rangle, ({\bf A}^\ell)_{ij}=\langle {\bf A}{\boldsymbol{\psi}}_i,{\boldsymbol{\psi}}_j\rangle \in\R^{\ell\times\ell}$ and ${\bf y_0^\ell}=({\bf {\bf \Psi}})^T{\bf y}_0\in\R^\ell$. We also note that ${\bf M}^\ell, {\bf A}^\ell\in\R^{\ell\times\ell}$.
The system (\ref{pod_sys}) is achieved following a Galerkin projection. If the dimension of the system is $\ell\ll n$, then a significant dimensionality reduction is accomplished.\\
This section focuses on several model order reduction techniques as they constitute the building blocks of the proposed method. In particular, we recall three key innovations for model reduction:  POD, DEIM and DMD.  These techniques provide an efficient projector for the reduction of the complexity of the problem under consideration.

%%%%%%%%%%%%%%%%%%%%%%%%%%%%%%%%%%%%%%%%%%%%%%%%%%%%%
\subsection{The POD method and reduced-order modeling}
\label{Section2.1}
%%%%%%%%%%%%%%%%%%%%%%%%%%%%%%%%%%%%%%%%%%%%%%%%%%%%%
One popular method for reducing the complexity of the system is the so-called Proper Orthogonal Decomposition (POD). The idea was proposed by Sirovich~\cite{Sir87} and is detailed here for completeness. We build an equidistant grid in time with constant step size $\Delta t$. Let $t_0:=0<t_1<t_2<\ldots<t_m\leq T$ with $t_j=j\Delta t,\;\;j=0,\ldots,m$. Let us assume we know the exact solution of (\ref{ode}) on the time grid points $t_j$, $j\in \{0,\ldots,m\}$. Our aim is to determine a POD basis of rank $\ell\ll n$ to optimally describe the set of data collected in time by solving the following minimization problem:
\begin{equation}\label{pbmin}
\min_{ {\boldsymbol{\psi}}_1,\ldots,{\boldsymbol{\psi}}_\ell\in\R^n} \sum_{j=0}^m \alpha_j\left\|{\bf y}(t_j)-\sum_{i=1}^\ell \langle {\bf y}(t_j),{\boldsymbol{\psi}}_i\rangle{\boldsymbol{\psi}}_i\right\|^2\quad \mbox{such that }\langle {\boldsymbol{\psi}}_i,{\boldsymbol{\psi}}_j\rangle=\delta_{ij},
\end{equation}
where the coefficients $\alpha_j$ are non-negative and ${\bf y}(t_j)$ are the so called {\em snapshots}, e.g. the solution of \eqref{ode} at a given time $t_j$.  Additionally, we assume ${\bf y}(t_j)\in V$ for a suitable Hilbert space $V$. The norm, here and in the sequel of the section, can be interpreted as the weighted norm such that $\langle {\bf u}, {\bf v}\rangle={\bf u}^T{\bf v}$ and $\|\cdot\|^2=\langle\cdot,\cdot\rangle$.

Solving (\ref{pbmin}) we look for  an orthonormal basis $\{{\boldsymbol{\psi}}_i\}_{i=1}^\ell$ which minimizes the distance between the sequence ${\bf y}(t_j)$ with respect to its projection onto this unknown basis. The matrix ${\bf Y}$ contains the collection of snapshots ${\bf y}(t_j)$ as columns. It is useful to look for $\ell\ll\min\{m,n\}$ in order to reduce the dimension of the problem considered.
The solution of (\ref{pbmin}) is given by the singular value decomposition of the snapshots matrix ${\bf Y}={\bf{\bf \Psi}}{\bf\Sigma} {\bf V}^T$, where we consider the first $\ell-$ columns $\{{\boldsymbol{\psi}}_i\}_{i=1}^\ell,$ of the orthogonal matrix ${\bf {\bf \Psi}}$ and set  ${\bf {\bf \Psi}}^{\mbox{\tiny POD}}={\bf \Psi}$.
%\end{theorem}

To concretely apply the POD method, the choice of the truncation parameter $\ell$ plays a critical role. There are no a-priori estimates which guarantee the ability to build a coherent reduced model, but one can focus on heuristic considerations, introduced by Sirovich \cite{Sir87}, so as to have the following ratio close to one:
\begin{equation}\label{indicator}
\mathcal{E}(\ell)=\dfrac{\sum\limits_{i=1}^\ell\sigma^2_i}{\sum\limits_{i=1}^r\sigma^2_i}.
\end{equation}
This indicator is motivated by the fact that the error in \eqref{pbmin} is given by the singular values we neglect:
\begin{equation}\label{err_pod}
\sum_{j=1}^m \alpha_j\left\|{\bf y}(t_j)-\sum_{i=1}^\ell \langle {\bf y}(t_j),{\boldsymbol{\psi}}_i\rangle{\boldsymbol{\psi}}_i\right\|^2=\sum_{i=\ell+1}^r \sigma_i^2,
\end{equation}
where $r$ is the rank of the snapshot matrix ${\bf Y}$. We note that the error \eqref{err_pod} is strictly related to the computation of the snapshots and it is not related to the reduced dynamical system.
More recently, Gavish and Donoho~\cite{gavish} have introduced a hard-thresholding technique for 
determining the truncation of the SVD when the data contains a low-rank signal with noise.  This method provides
a principled approach to rank selection.

%Let us assume the POD basis is computed, one may obtain a reduced order model for the problem \eqref{ode}. 

%%%%%%%%%%%%%%%%%%%%%%%%%%%%%%%%%%%%%%%%%%%%%%%%%%%%%
\subsection{Discrete Empirical Interpolation Method}
\label{Section2.2}
%%%%%%%%%%%%%%%%%%%%%%%%%%%%%%%%%%%%%%%%%%%%%%%%%%%%%

%%%%%%%%%%%%%%%%%%%%%%%%%%%%%%%%%%%%%%%%%%%%%%%%%%%%%

For a review of DEIM, we closely follow the presentation in \cite{Vol11}. The ROM introduced in (\ref{pod_sys}) is a nonlinear system where the significant challenge with the POD--Galerkin approach is the computational complexity associated with the evaluation of the nonlinearity. To illustrate this issue, we consider the nonlinearity in (\ref{pod_sys}): 
$$ {\bf F}(t,{\bf y}^\ell(t))=({\bf {\bf \Psi}}^{\mbox{\tiny POD}})^T {\bf f}(t,{\bf {\bf \Psi}}^{\mbox{\tiny POD}} {\bf y}^\ell(t))=\langle {\bf f}(t,{	\bf y}(t)), {\bf {\bf \Psi}}^{\mbox{\tiny POD}}\rangle.$$
To compute this inner product, the variable ${\bf y}^\ell(t)\in\R^\ell$ is first expanded to an $n-$dimensional vector ${\bf {\bf \Psi}}^{\mbox{\tiny POD}} {\bf y}^\ell(t)\in\R^n$, then the nonlinearity ${\bf f}(t,{\bf {\bf \Psi}}^{\mbox{\tiny POD}} {\bf y}^\ell(t))$ is evaluated and, at the end, we return back to the reduced-order model. This is computationally expensive since it implies that the evaluation of the nonlinear term requires computing the full, high-dimensional model, and therefore the reduced model is not independent of the full dimension $n.$ %We note that, for simplicity, we dropped the weighted inner product.
To avoid this computationally expensive, high-dimensional evaluation, the gappy POD method was introduced~\cite{gap1}.
In its original formulation, random and sparse sampling was proposed for computing the required nonlinear inner products.  Advances in gappy methods have led to the state-of-the-art
{\em Empirical Interpolation Method} (EIM, \cite{BMNP04}) and {\em Discrete Empirical Interpolation Method} (DEIM, \cite{CS10}) methods which are now broadly used in the ROMs community. %We note that DEIM is built upon EIM:  the two methods are essentially equivalent and are based on a POD approach combined with a greedy algorithm. DEIM is the tensorial matricial version of EIM, and it is used here due to the nature of our time-dependent problem. The interested reader is referred to \cite{CS10} for further information. 

The computation of the POD basis functions for the nonlinear part are related to the set of the snapshots ${\bf f}(t_j,{\bf y}(t_j))$ where ${\bf y}(t_j)$ is already computed from \eqref{ode}. We denote with ${\bf U}\in\R^{n\times k}$ the POD basis function of rank $k$ of the nonlinear part. The DEIM approximation of ${\bf  f}(t,{\bf y}(t))$ is as follows
$${\bf f}^{\mbox{\tiny DEIM}}(t,{\bf y}^{\mbox{\tiny DEIM}}(t))={\bf U}({\bf S}^T {\bf U})^{-1} {\bf f}(t,{\bf y}^{\mbox{\tiny DEIM}}(t))$$
where ${\bf S}\in\R^{n\times k}$ and ${\bf y}^{\mbox{\tiny DEIM}}(t)={\bf S}^T{\bf {\bf \Psi}}^{\mbox{\tiny POD}}{\bf y}^\ell(t)$. The matrix ${\bf S}$ is the interpolation point where the nonlinearity is evaluated and the selection of its points is made according to an LU decomposition algorithm with pivoting~\cite{CS10}, or following the QR decomposition with pivoting~\cite{DG15}. 
The error between ${\bf f}(t,{\bf y}(t))$ and its DEIM approximation $f^{\mbox{\tiny DEIM}}$ is given by
$$\|{\bf f}-{\bf f}^{\mbox{\tiny DEIM}}\|_2\leq  c\|({\bf I}-{\bf UU}^{T})f\|_2\quad \,\,\, \mbox{with} \,\,\, c=\|({\bf S}^T{\bf U})^{-1}\|_2$$
where different error performance is achieved depending on the selection of the interpolation points in $S$ as shown in \cite{DG15}.

\subsection{Dynamic Mode Decomposition}
\label{Section2.3}
%\setcounter{section}{3}
%\setcounter{equation}{0}
%\setcounter{theorem}{0}
%\setcounter{algorithm}{0}
%\renewcommand{\theequation}{\arabic{section}.\arabic{equation}}
%%%%%%%%%%%%%%%%%%%%%%%%%%%%%%%%%%%%%%%%%%%%%%%%%%%%%

DMD is an {\em equation-free}, data-driven method capable of providing accurate assessments of the spatio-temporal coherent structures in a given complex system, or short-time future estimates of such a systems. It traces its origins to pioneering work of Bernard Koopman in 1931~\cite{koopman}, whose work was revived in a set of papers starting in 2004~\cite{Mezic2004,Mezic2005,mezic2}. The DMD provides the eigenvalues and eigenvectors of the best fit linear system relating a snapshot matrix and a time shifter version of the snapshot matrix at some later time.

Consider the following data snapshot matrices
%
%\begin{subeqnarray}
\begin{equation}
  {\bf Y} \!=\! \begin{bmatrix}
\vline & \vline & & \vline \\
{\bf y}(t_0) & {\bf y}(t_1) & \cdots & {\bf y}(t_{m-1})\\
\vline & \vline & & \vline
\end{bmatrix}, \hspace{0.1in} {\bf Y}' \!=\! \begin{bmatrix}
\vline & \vline & & \vline \\
{\bf y}(t_1) & {\bf y}(t_2) & \cdots & {\bf y}(t_m)\\
\vline & \vline & & \vline
\end{bmatrix}
\end{equation}
%\end{subeqnarray}
%
with ${\bf y}(t_j)$ an initial condition to \eqref{ode} and ${\bf y}(t_{j+1})$ its corresponding output
after some prescribed evolution time $\Delta t>0$ with there being $m$ initial conditions considered. The DMD involves the decomposition of the best-fit linear operator ${\bf A}$ relating the matrices above:
\begin{equation}
{\bf Y}'={\bf AY}.
\end{equation}
where ${\bf A}\in\R^{n\times n}$ is unknown. The exact DMD algorithm proceeds as follows \cite{DMD5}:
First, we collect data {\bf Y, Y'} and compute the reduced singular value decomposition of {\bf Y}:
$${\bf Y}={\bf U}{\bf \Sigma}{\bf V}^*.$$ Then, compute the least-squares fit ${\bf A}$ that satisfies ${\bf Y}'={\bf AY}$ and project onto POD modes ${\bf U:}$
$${\bf \tilde{A}}={\bf U}^* {\bf AU}={\bf U}^*{\bf Y}'{\bf V}{\bf \Sigma}^{-1},$$
and compute the eigen-decomposition of ${\bf \tilde{A} }:$
$${\bf\tilde A}{\bf W}={\bf W \Lambda},$$
where ${\bf \Lambda}$ are the DMD eigenvalues.
The DMD modes ${\bf \Psi}^{\mbox{\tiny DMD}}$ are given by:
\begin{equation}\label{dmd_basis}
{\bf\Psi}^{\mbox{\tiny DMD}}={\bf Y}'{\bf V}{\bf \Sigma}^{-1}{\bf W}.
\end{equation}
The data ${\bf Y, Y}'$ may come from a nonlinear system $y(t_{j+1})=f(y(t_j)),$
in which case the DMD modes are related to eigenvectors of the infinite-dimensional Koopman operator. More details can be found in~\cite{dmdbook}. We may interpret DMD as a model reduction technique if data is acquired from a high-dimensional model, or a method of system identification if the data comes from measurements of an unknown system. For the purpose of this work we consider the DMD-Galerkin method, where the assumption \eqref{pod_ans} hold true for ${\bf \Psi}$ given by \eqref{dmd_basis}. We note that the techniques we provide aim to speed up the computation of the offline stage, whereas the online stage will present the same cost as the standard methods.

\section{Randomized Linear Algebra in Model Order Reduction}\label{Section3}

Randomized linear algebra is of growing importance for the analysis of high-dimensional data~\cite{gunnar_rLA}. Specifically, randomized techniques attempt to construct low-rank matrix decompositions that are computationally efficient and accurate approximations of the standard matrix decompositions such as QR and SVD.   Randomized algorithms can be parallelized and distributed for large matrices and there are several implementations of the randomized techniques in \textit{MATLAB} or \textit{R} that are now available via open source \cite{erichson2016,szlam,voronin2015}.  
 The algorithms
that result from using randomized sampling techniques are not only computationally efficient, but are also
simple to implement as they rely on standard matrix-matrix multiplication and unpivoted QR factorization.

\begin{figure}[t]
\begin{center}
\vspace*{-.2in}
\includegraphics[scale=0.6]{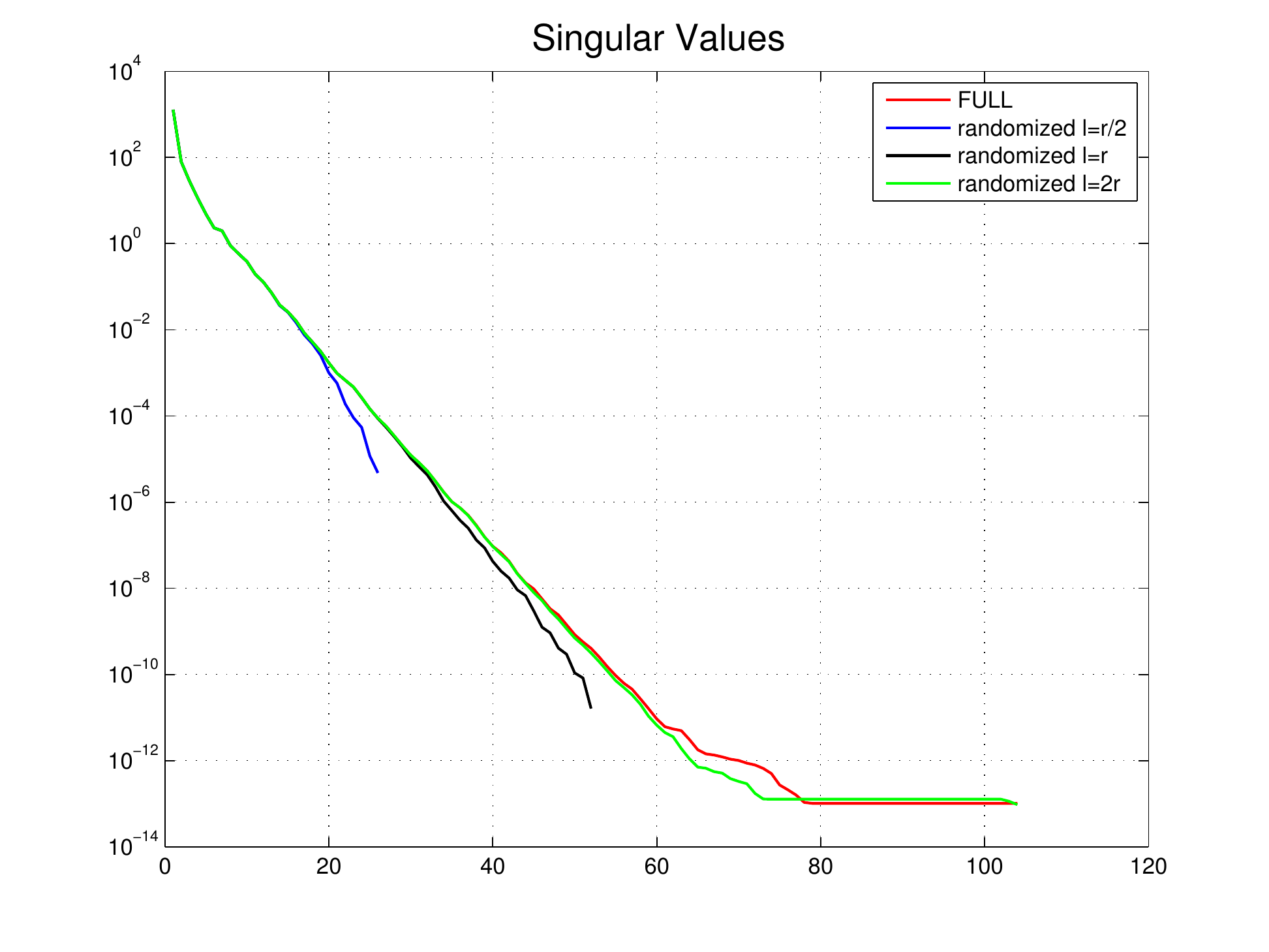}
\end{center}
\vspace*{-.2in}
\caption{Convergence of the singular values for randomized SVD with different number of measurements $\ell$. We consider the full matrix (red), $\ell=r/2$ (blue), $\ell=r$ (black), $\ell=2r$ (green) where $r$ is the rank of matrix.}
\label{fig:rsvd}
\end{figure}

Consider a randomized algorithm to compute the low-rank matrix approximation~\cite{halko2011rand}
\begin{equation*}\label{eq:QB}
\begin{array}{cccc}
\mathbf{A} & \approx & \mathbf{Q} & \mathbf{B} \\
n\times m &   &  n\times \ell & \ell\times m
\end{array} 
\end{equation*}
where $\ell$ denotes the target-rank and is assumed to be $\ell\ll \min\{m,n\}$.  Random matrix theory provides a simple and elegant solution for computing the low-rank approximation by creating a random sampling matrix $\mathbf{\Omega} \in \mathbb{R}^{n\times \ell}$ where the entries are drawn from, for example, a Gaussian distribution. Then, a sampled matrix $\mathbf{Y} \in \mathbb{R}^{n\times \ell}$ is computed as
\begin{equation*}\label{eq:YAQ}
\mathbf{Y} = \mathbf{A}\mathbf{\Omega}.
\end{equation*}
If the matrix $\mathbf{A}$ has exact rank $\ell$, then the sampled matrix $\mathbf{Y}$ spans, with high probability, a basis for the column space.  However, most data matrices in practice are only dominated by rank-$\ell$ features since the singular values $\{\sigma_i\}_{i=\ell+1}^n$ are non-zero. Thus, instead of just using $\ell$ samples, it is favorable to slightly oversample $\ell=\ell+p$, were $p$ denotes the number of additional samples. In practice, small values of $p\approx 5-10$ are sufficient to obtain a good basis that is comparable to the best possible basis~\cite{gunnar_rLA}. An orthonormal basis $\mathbf{Q} \in \mathbb{R}^{n\times \ell}$ is then obtained via the QR-decomposition  $\mathbf{Y}=\mathbf{Q}\mathbf{R}$, such that
\begin{equation*}
\mathbf{A} \approx \mathbf{Q}\mathbf{Q}^\intercal\mathbf{A} .
\end{equation*}
Finally, $\mathbf{A}$ is projected to this low-dimensional space
\begin{equation*}
\mathbf{B} = \mathbf{Q}^\intercal\mathbf{A},
\end{equation*}
where $\mathbf{B} \in \mathbb{R}^{\ell \times m}$. 
The matrix $\mathbf{B}$ can then be used to efficiently compute the matrix decomposition of interest such as the SVD. The oversampling $p$ allows one to control the approximation error\cite{halko2011rand,gunnar_rLA}.  The algorithm is summarized in Algorithm \ref{alg:rSVD}.  In Figure \ref{fig:rsvd} we show the decay of the singular values for different level of the randomized SVD. As expected increasing the number of sampling we obtain more accurate approximation.

\begin{algorithm}[t]
\caption{Randomized SVD (rSVD)}
\label{alg:rSVD}
\begin{algorithmic}[1]
\REQUIRE Matrix ${\bf Y}\in\R^{n\times m}$
\STATE Draw a Gaussian random matrix ${\bf \Omega}\in\R^{m\times \ell}$
\STATE Form the sample matrix ${\bf X}= {\bf \Omega Y}$
\STATE Compute the QR decomposition of ${\bf X}: {\bf X}={\bf QR}$
\STATE Set {\bf B}=${\bf Q}^T{\bf Y}$
\STATE Compute SVD of ${\bf B}={\bf \hat U}{\bf \Sigma} {\bf V}^T$
\STATE Set ${\bf U}={\bf Q \hat U}$
\end{algorithmic}
\end{algorithm}

\subsection{Compressed Model Order Reduction Techniques }
\label{Section4}
\setcounter{equation}{0}
\setcounter{theorem}{0}
\setcounter{algorithm}{0}
\renewcommand{\theequation}{\arabic{section}.\arabic{equation}}
%%%%%%%%%%%%%%%%%%%%%%%%%%%%%%%%%%%%%%%%%%%%%%%%%%%%%
Model order reduction techniques are usually based on snapshots that collect data on the underlying dynamical system. The SVD decomposition of the date matrix ${\bf Y}\in\R^{n\times m}$ provides a low-dimension projector operator that allows one to obtain surrogate models. However, the SVD may be computationally expensive and, for this reason, we propose the use of Algorithm \ref{alg:rSVD} to reduced the offline cost of the method.
The main idea is to consider basis functions not from the full set of measurements but from a few spatially incoherent measurements. We introduce the measurement matrix ${\bf C}\in\R^{p\times n}$ which produces the compressed matrix ${\bf X}\in\R^{p\times m}$ such that:
$${\bf X}={\bf C Y}.$$
Here, we consider sparse measurements of the snapshots matrix in order compute POD and DMD from this new compressed snapshot matrix. In this paper we assume that snapshots matrix is almost square, e.g. $n\approx m$, and one can imagine this is a realistic situation working with an explicit time scheme or in a many-query context. In the following subsections we provide further details about compressed POD, compressed POD-DEIM, and compressed DMD.

\subsubsection{Compressed POD}
The compressed POD method works, as POD, starting with a snapshot matrix with the aim to compute solutions of the problem \eqref{pbmin} in a fast and reliable way. As discussed before, the solution of the minimization problem leads to an expensive singular value decomposition problem.  Here, the idea is to apply the randomized SVD technique in Section 2 for the approximation of \eqref{pbmin}. The method works as follows: (i) we collect the snapshot set and (ii) we solve the optimization problem \eqref{pbmin}.   We make use of the optimality conditions in \cite{Vol11} in order to take advantage of the randomized SVD. In this way we are able to compute the compressed POD basis functions in a significantly faster way. Clearly the number of samples point plays a crucial role. The algorithm is summarized in Algorithm \ref{alg:cPOD}.

\begin{algorithm}
\caption{Compressed POD (cPOD)}
\label{alg:cPOD}
\begin{algorithmic}[1]
\REQUIRE Snapshot Matrix ${\bf Y}\in\R^{n\times m}$, $\ell$ number of basis functions., $p$ number of measurements.
\STATE Compute the Randomized SVD (see Algorithm \ref{alg:rSVD}), $[{\bf U}, {\bf \Sigma}, {\bf V}]=rsvd({\bf Y})$
\STATE Set ${\bf \Psi}_i= {\bf U}_i$ for $i=1,\ldots,\ell$.
\end{algorithmic}
\end{algorithm}

The error in the minimization problem is now associated with subsampling of  the randomized SVD (\cite{halko2011rand}):
\begin{equation}\label{err_cpod}
\mathbb{E}\left(\sum_{j=1}^m \alpha_j\left\|{\bf y}(t_j)-\sum_{i=1}^\ell \langle {\bf y}(t_j),{\boldsymbol{\psi}}_i\rangle{\boldsymbol{\psi}}_i\right\|^2\right)=\left(1+\sqrt{\dfrac{\ell}{p-1}}\right)\sigma_{\ell+1}^2+\dfrac{\sqrt{\ell+p}}{p}\sum_{j=\ell+1}^d\sigma_j^2.
\end{equation}
where $d$ is the rank of the snapshot matrix ${\bf Y}$ and $p$ is the number of the samples in the compressed technique. We note that we consider the expectation value of the error due to the random measurements we consider. The error \eqref{err_cpod} is now related to the computation of the set of snapshots and the number of samples $p$. We note that if the singular values of the snapshot matrix decay rapidly a minimal amount of samples drives the error close to the theoretically minimum value.  However, if the singular values do not decay rapidly we can lose accuracy. We refer to \cite{halko2011rand} and the reference therein for more details about the error of the randomized SVD. Finally, we note that the POD basis functions for the snapshot matrix ${\bf Y}\in\R^{n\times m}$ can be also computed from the eigenvalue problem for the matrix ${\bf Y} {\bf Y}^T$ or ${\bf Y}^T{\bf Y}$.  Interested readers can see Ref.~\cite{Vol11} for a more comprehensive description.  Regardless, if the dimension of the matrix ${\bf Y}$ is such that $n\approx m$, then the computation of the eigenvalue problem will not lead to a faster approximation than the SVD. This further motivates our approach through the rSVD.

\subsubsection{Compressed POD-DEIM}
Similarly to the compressed POD method we aim to apply the rSVD to the DEIM approach. The DEIM method considers the computation of the SVD for both the snapshots of the solution and snapshots of the nonlinear term. We note that, although the online stage benefits from a sparse evaluation of the nonlinearity, the offline stage is even more expensive than POD itself. The goal is to substitute the full dimensional SVD with the much smaller randomized SVD. In this way we can highly reduce the cost of the computational costs and, at the same time, obtain accurate results.

%%%%%%%%%%%%%%%%%%%%%%%%%%%%%%%%%%%%%%%%%%%%%%%%%%%%%
%%%%%%  SECTION: Dynami Mode Decomposition
%%%%%%%%%%%%%%%%%%%%%%%%%%%%%%%%%%%%%%%%%%%%%%%%%%%%%

\subsection{Compressive DMD}
We can also combine ideas from compressive sampling to compute the dynamic mode decomposition from a few measurements of the data. This method was already introduced in \cite{brunton:cs}. Here, it is applied in as Galerkin projection method. It is possible to either collect data ${\bf Y, Y}'$ or projected data ${\bf X, X}',$ where ${\bf X}={\bf CY}$, ${\bf X}'={\bf CY}'$ and ${\bf C}\in\R^{p\times m}$ is the measurement matrix. We will call the matrices ${\bf X, X}'$ the output-projected snapshot matrices.
Similar to equation above, ${\bf X}$ and ${\bf X}'$ are related by
\begin{equation}
{\bf X}'={\bf A_X X}.
\end{equation}
The goal, as in DMD, is to compute eigenvalues and eigenvectors of the unknown matrix ${\bf A_X}$. The method differs from the standard DMD since we are using sparse measurements. Under general assumptions it is possible to prove the convergence of the method when the number of measurements $p$ increases. Cleary the cDMD method is computationally more efficient, and the method is summarized in Algorithm \ref{Alg_cDMD}.

\begin{algorithm}
\caption{Compressive DMD (cDMD)}
\label{Alg_cDMD}
\begin{algorithmic}[1]
\REQUIRE Snapshots $\{{\bf y}(t_0),\ldots,{\bf y}(t_m)\}$, ${\bf C}\in\R^{p\times m}$
\STATE Set ${\bf Y}=[{\bf y}(t_0),\dots, {\bf y}(t_{m-1})]$ and $Y'=[{\bf y}(t_1),\dots, {\bf y}(t_m)]$,
\STATE ${\bf X}={\bf CY}$, ${\bf X}'={\bf CY}'$
\STATE Compute the SVD of ${\bf X}$, ${\bf X}={\bf U}{\bf \Sigma}{\bf  V}^T$
\STATE Define $\tilde{{\bf A}}_{\bf x}:={\bf U}^*{\bf Y}'{\bf V}{\bf \Sigma}^{-1}$
\STATE Compute eigenvalues and eigenvectors of $\tilde{{\bf A}}_{\bf x} {\bf W}={\bf W}{\bf \Lambda}$.
\STATE Set ${\bf \Psi}^{\mbox{\tiny DMD}}={\bf X}'{\bf V}{\bf \Sigma}^{-1}{\bf W}$ \\
\end{algorithmic}
\end{algorithm}
Once the DMD basis functions ${\bf \Psi}^{\mbox{\tiny DMD}}$ are computed we utilize assumption \eqref{pod_ans} and obtained a surrogate model of the form \eqref{pod_sys}.

%%%%%%%%%%%%%%%%%%%%%%%%%%%%%%%%%%%%%%%%%%%%%%%%%%%%%
\section{Numerical Tests}
\label{Section5}
\setcounter{section}{5}
\setcounter{equation}{0}
\setcounter{theorem}{0}
\setcounter{algorithm}{0}
\renewcommand{\theequation}{\arabic{section}.\arabic{equation}}
%%%%%%%%%%%%%%%%%%%%%%%%%%%%%%%%%%%%%%%%%%%%%%%%%%%%%

In this section we present our numerical tests using our three proposed compressed/randomized SVD strategies of the last section. In our numerical computations we use the finite difference method to reduce a partial differential equation into the form \eqref{ode} and integrate the system with a semi-implicit scheme in the first example and Newton method in the second. All the numerical simulations reported in this paper are performed on a MacBook Pro with an Intel Core i5, 2.2Ghz and 8GB RAM using MATLAB R2013a. 

In the following numerical examples we build different surrogate models, such as POD, compressed POD (cPOD), POD-DEIM,  compressed POD-DEIM (cPOD-cDEIM), DMD and compressed DMD (cDMD) and compare their performance in terms of CPU time and the error with respect to a reference solution computed by a high-fidelity, finite-difference approximation. We select two numerical examples, the first one considers a time-dependent semi-linear PDEs whereas the second studies a semi--linear elliptic parametric equations. Both examples lead to the same conclusions. In the numerical tests, the number samples utilized for the compression of the snapshot matrix is twice the rank. As shown in Figure \ref{fig:rsvd}, this turns out to be very efficient for both accuracy and computational cost.

\subsection{Test 1: Semi-Linear Equation}
\begin{figure}[htbp]
\begin{center}
\includegraphics[scale=0.25]{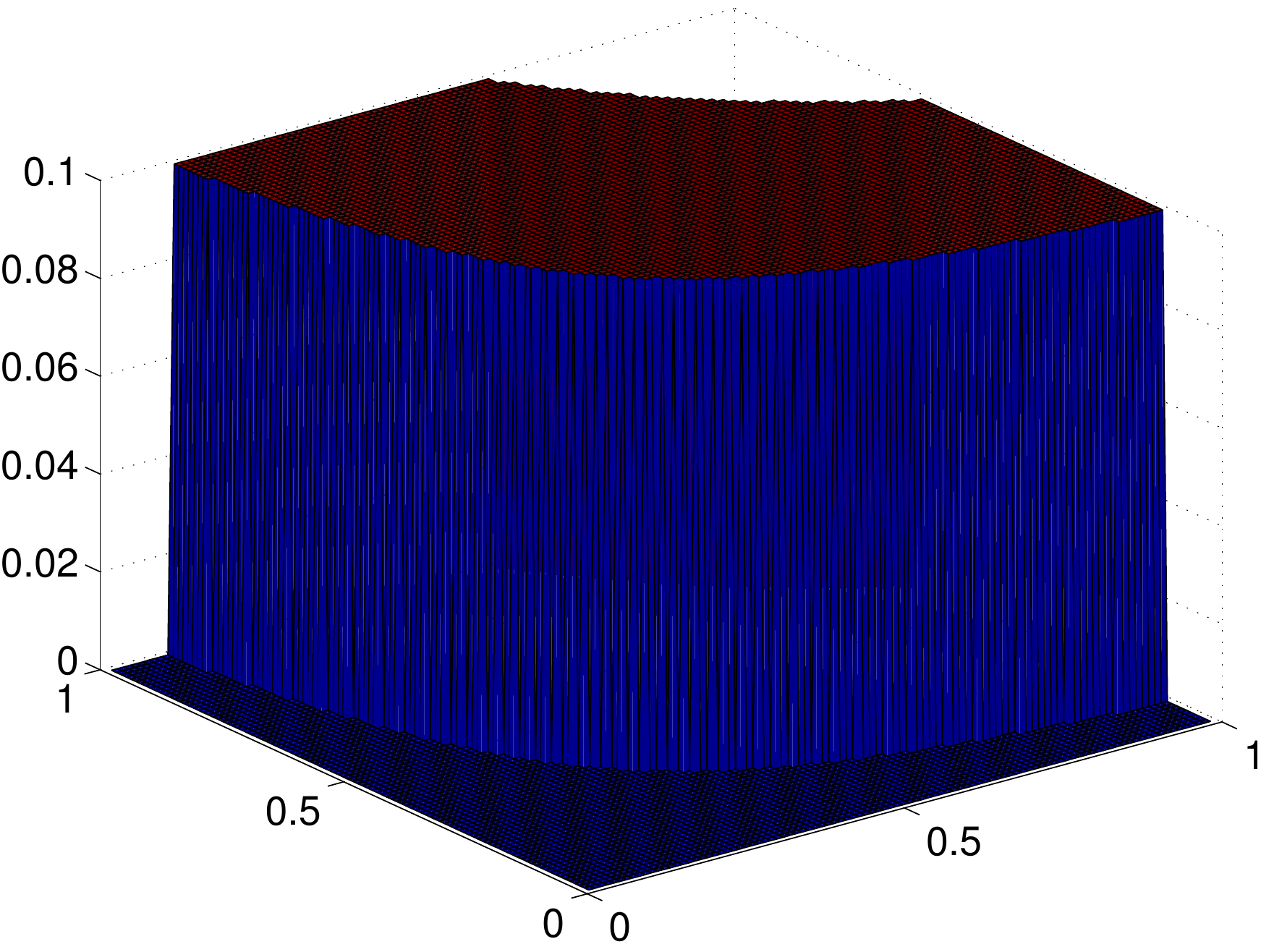}
\includegraphics[scale=0.25]{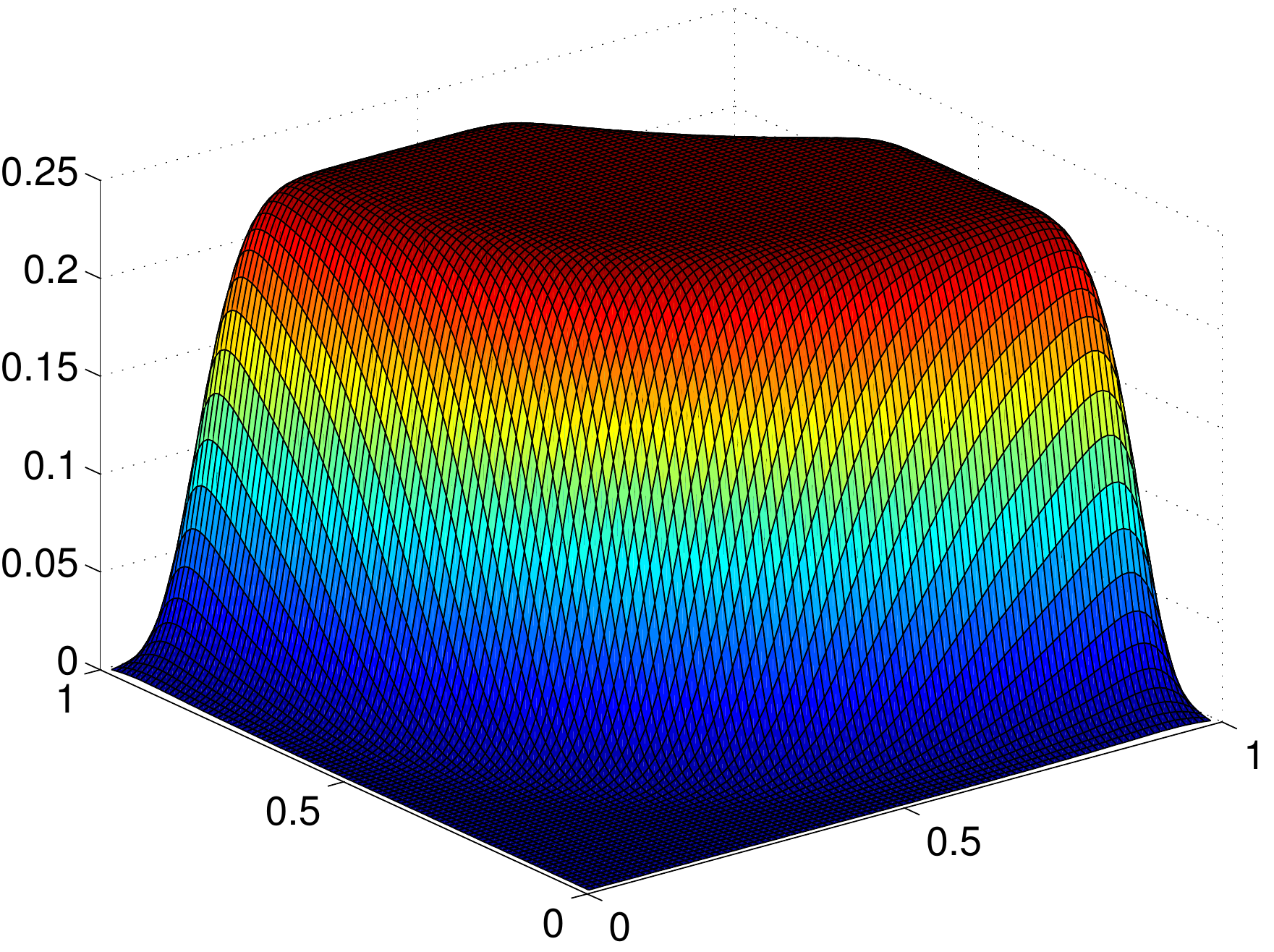}\\
\includegraphics[scale=0.25]{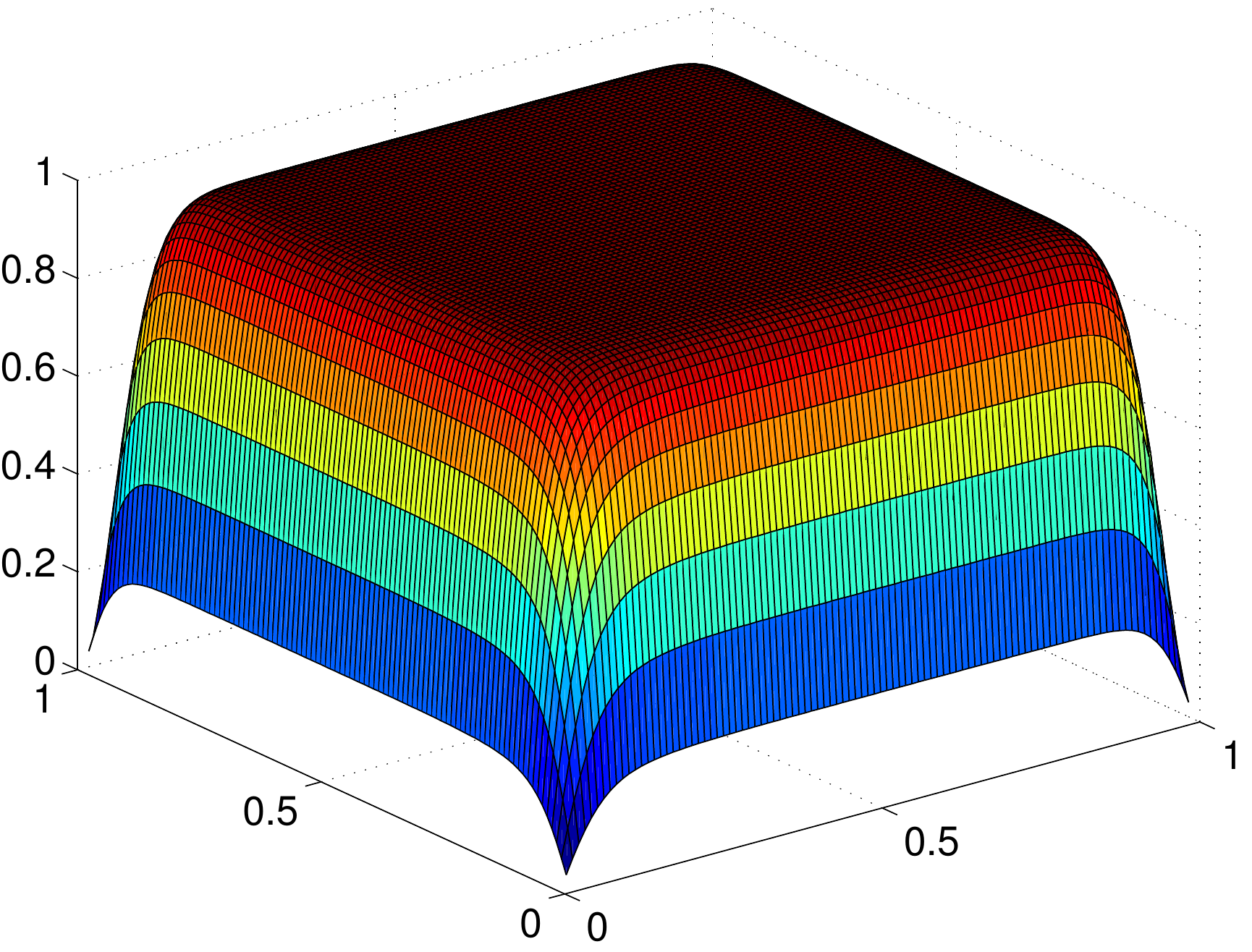}
\includegraphics[scale=0.25]{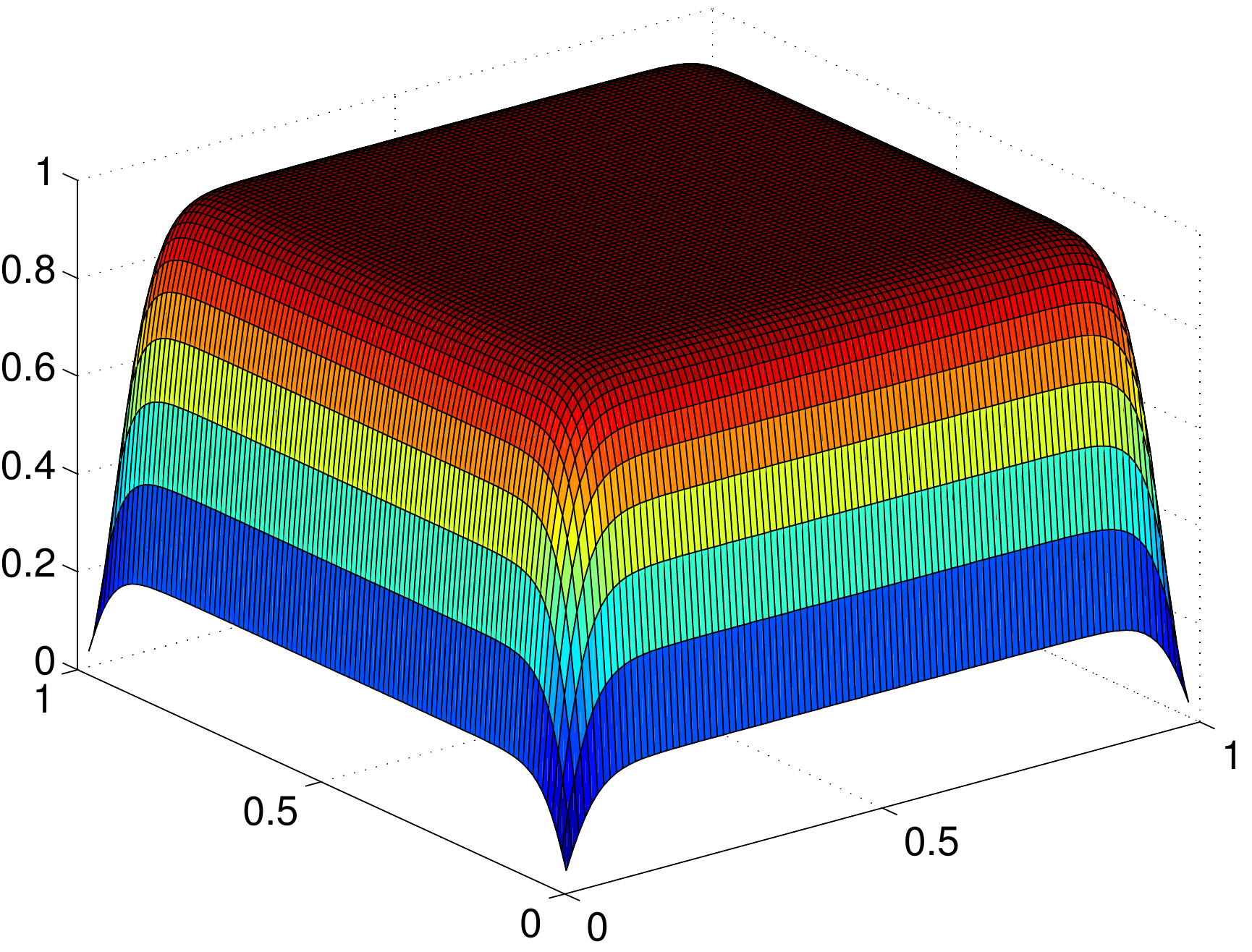}\\
\end{center}
\label{test1:sol}
\caption{Test 1: Solution of equation \eqref{prb_test1} at time $t=\{0,0.1\}$ (top) and $t=\{1.5,3\}$ (bottom).}
\end{figure}

Let us consider the following semi linear parabolic equation:
\begin{equation}\label{prb_test1}
\left\{
\begin{aligned}
y_t(x,t)-\theta \Delta y(x,t)+\mu(y-y^3)&=0, &\qquad (x,t)\in\Omega\times[0,T],\\
y(x,0)&=y_0(x), &\qquad x\in \Omega,\\
y(\cdot,t)&=0, &\quad x\in\partial\Omega, t\in [0,T],
\end{aligned}\right.
\end{equation}
where $\Omega=[0,1]\times[0,1], T=5, \theta=0.1, \mu=1, x=(x_1,x_2), y_0(x)=0.1$ if $0.1\leq x_1x_2\leq0.6$ and $0$ elsewhere. The POD basis vectors are built upon 10000 equidistant snapshots. The FD discretization yields a system of ODEs of the same form as \eqref{ode} with $n=10000$. The solution of this equation generates a stationary solution $y(x,t)\equiv 1$ for large $t$ as shown in Figure \ref{test1:sol}.
%Figure \ref{test1:svd} shows a similar decay for the singular values of the snapshot set, and of the nonlinear term in \eqref{prb_test1}.

\begin{figure}[htbp]
\begin{center}
\includegraphics[scale=0.25]{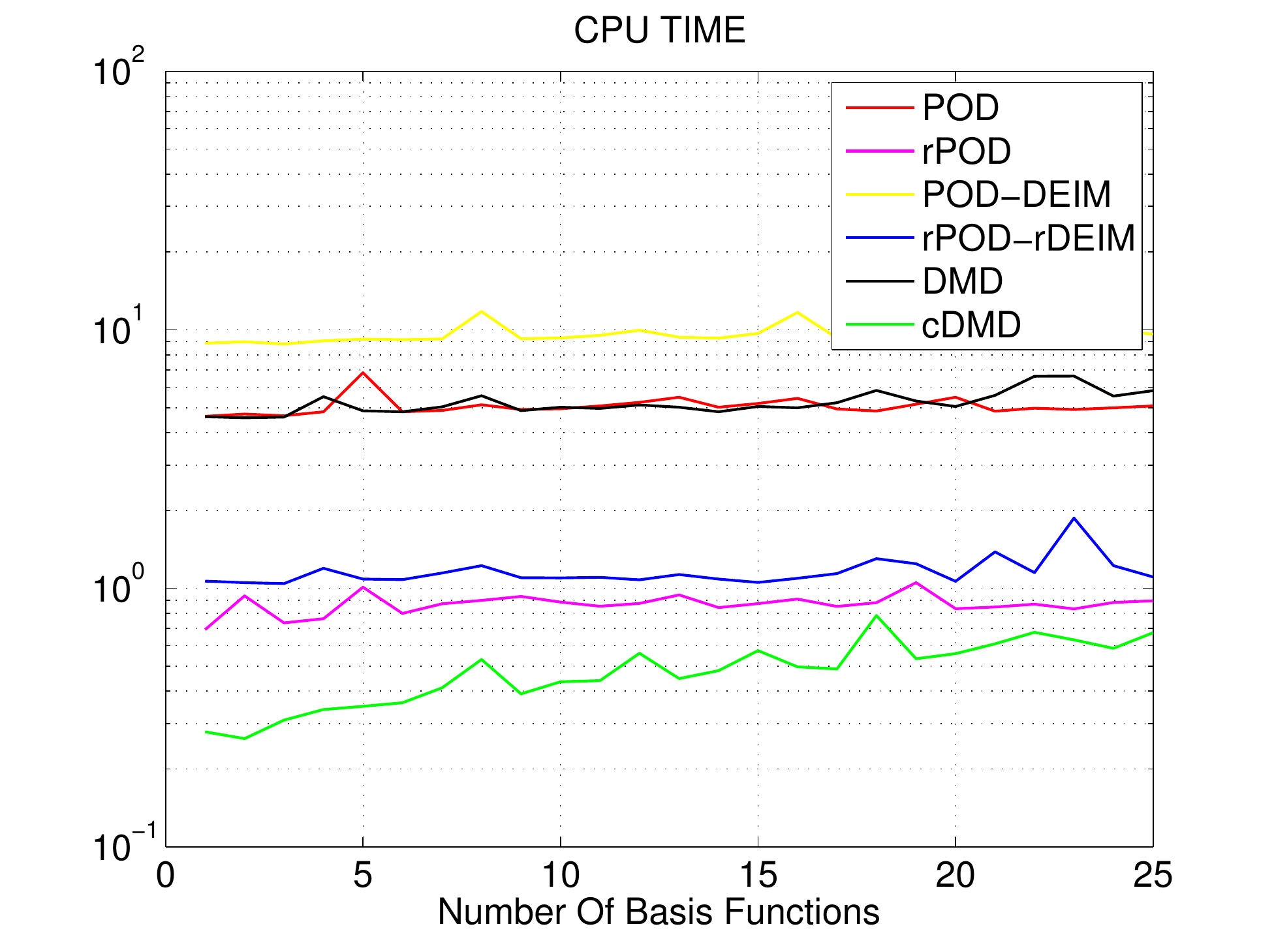}
\includegraphics[scale=0.25]{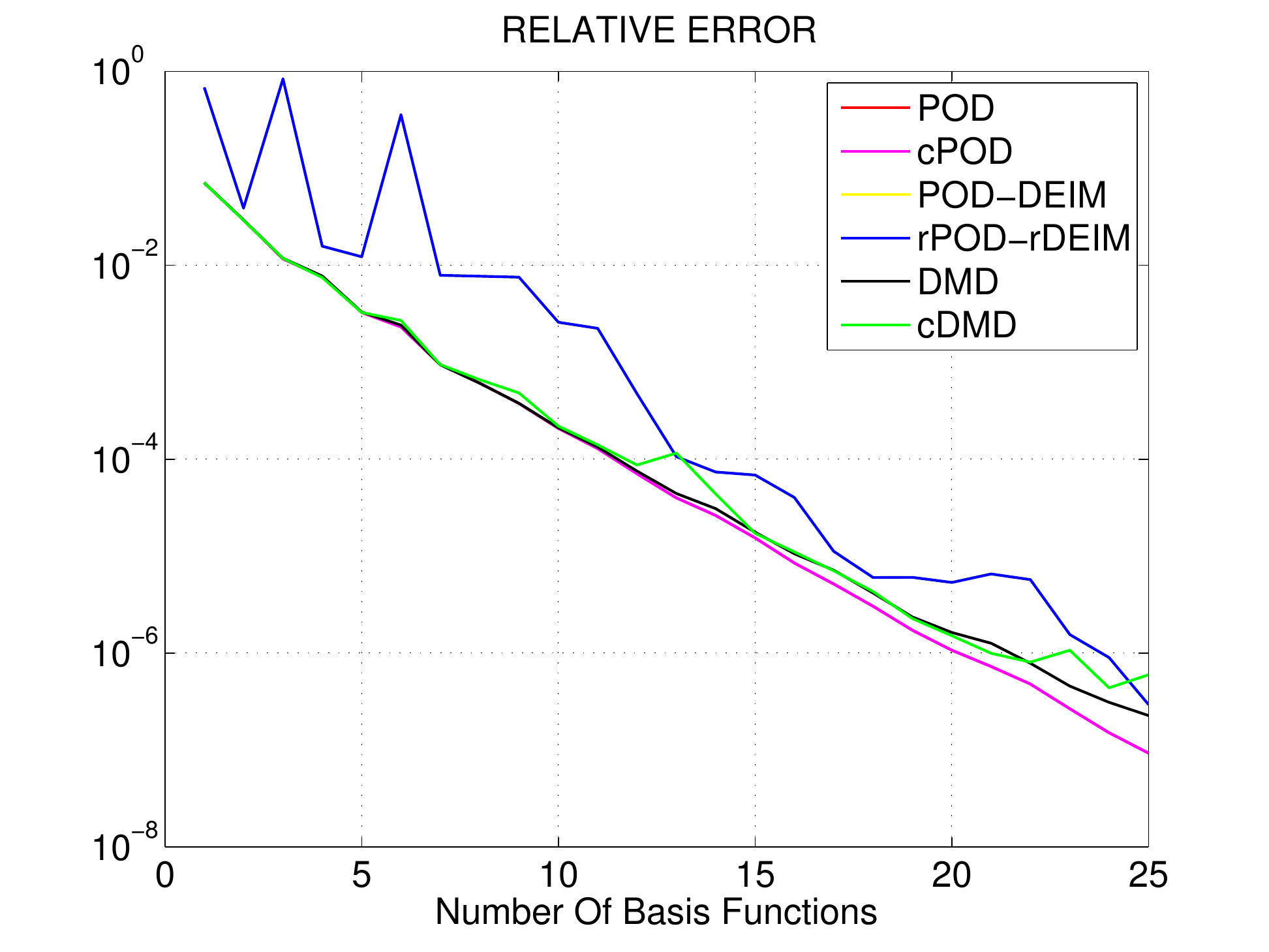}
\end{center}
\caption{Test 1: CPU-time of the offline-online stages (left) and Relative Error in Frobenius norm (right). We compare the following methods: POD (red), cPOD (magenta), POD-DEIM (yellow), cPOD-cDEIM (blue, DMD (black), cDMD (green). Number of model are always the same for all the methods.}
\label{test1:an}
\end{figure}

The complexity of problem \eqref{prb_test1} is reduced by model order reduction. When dealing with model order reduction, it is relevant to consider the CPU time of the simulation and the error. In general it is important to have a trade-off between the two quantities. Figure \ref{test1:an} considers the CPU time on the left panel. As we can see the compressed techniques are faster than the standard reduction techniques. We note here that for the CPU time we consider both offline and online stages. Although we do not aim at and improvement of the online stage, in this work one might also consider a further speed up as suggested in \cite{AK16}. It is somehow clear that the compressed DMD provide the fastest approximation since it does not require the computation of the randomized SVD. However, we show in Figure \ref{test1:an} the relative error computed with respect to the Frobenius norm. As we can see, POD and cPOD, such as POD-DEIM and cPOD-cDEIM, perform exactly the same results. Slightly different are the results from DMD and cDMD. However, all these techniques perform with very high accuracy. As expected, the POD-DEIM and its related compressed technique is less accurate since we do not evaluate the nonlinearity for the full state.

\begin{figure}[htbp]
\begin{center}
\includegraphics[scale=0.25]{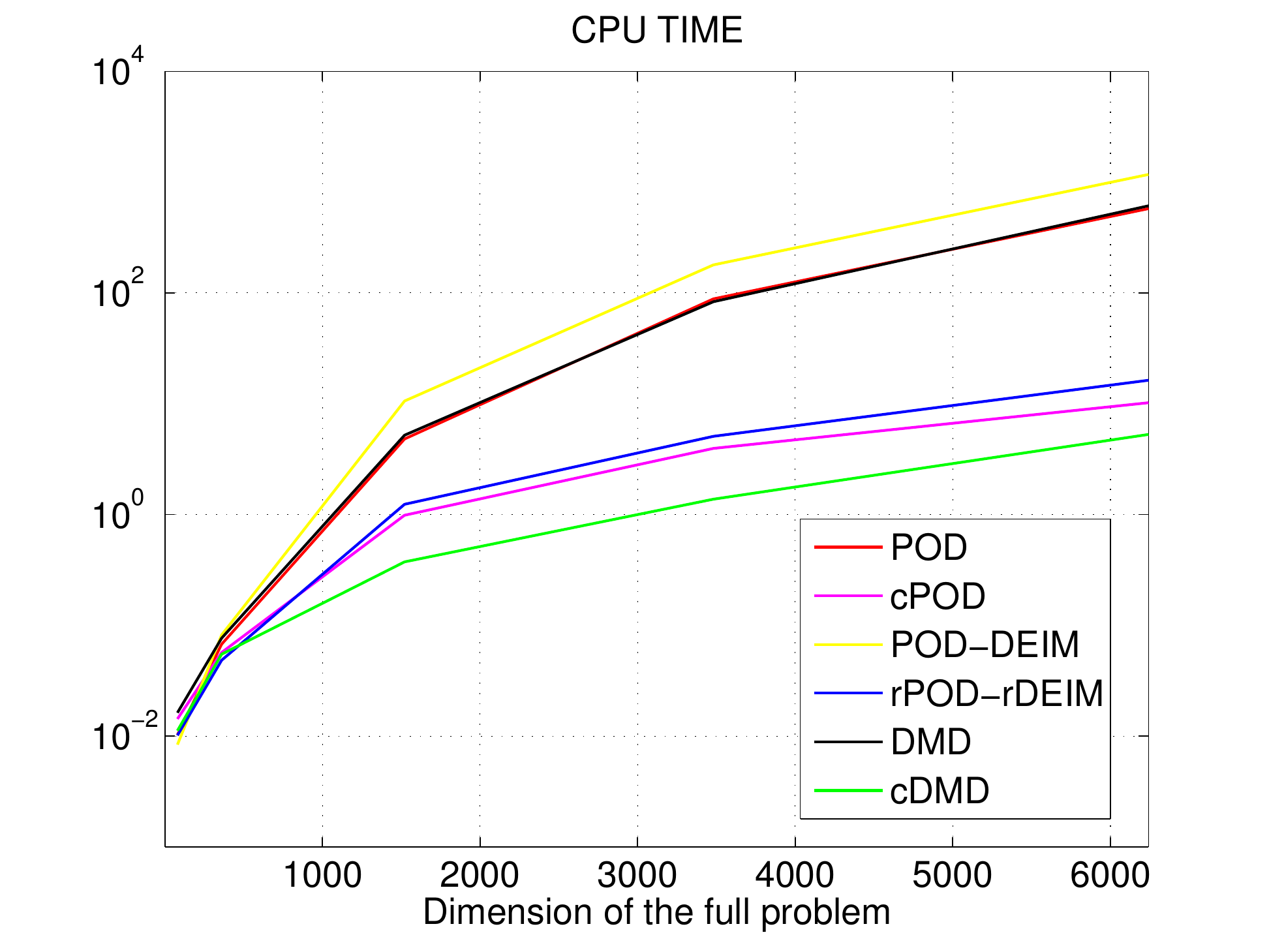}
\includegraphics[scale=0.25]{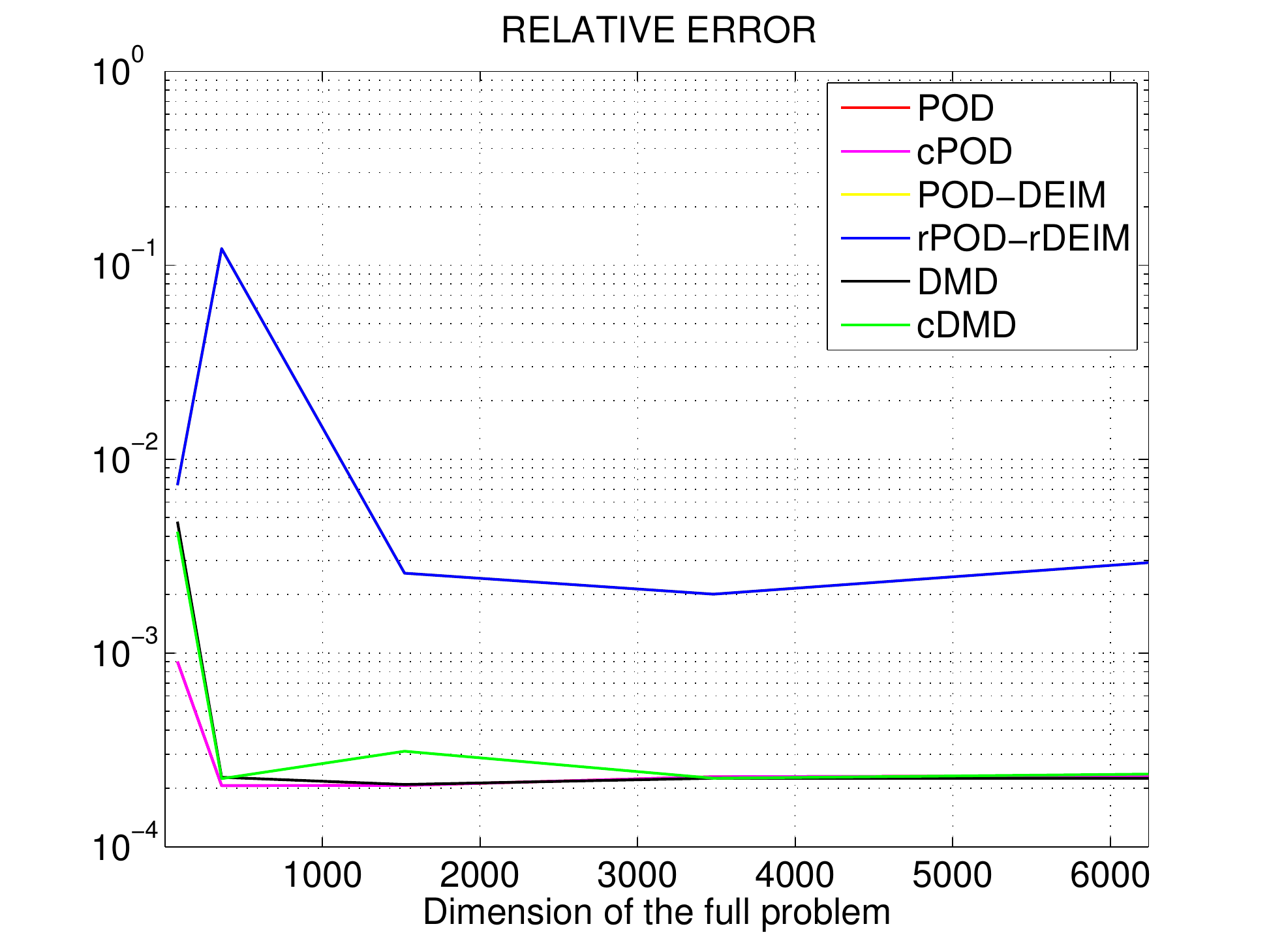}
\end{center}
\caption{Test 1: Scaling of the CPU time with increasing dimension of the snapshot set (left), Relative Error for 10 modes and different snapshot set.}
\label{test1:scale}
\end{figure}

Another important feature to investigate when dealing with compressed techniques is how the CPU time scales with different dimensions of the snapshot matrix. The computation of the SVD is, computationally, the most expensive part of the method and its cost varies according to the dimension of the snapshot set. Here we consider a square matrix. As we can see in Figure \ref{test1:scale} in the left panel, the CPU time scale shows that we gain more than 2 order of magnitude in speed up as the dimension increases. Thus, it provides a powerful technique that allows one to significantly reduce the computational costs in the offline stage. In the right panel we can see the relative error for 10 basis functions. 

\subsection{Test 2: Parametric Example}
\begin{figure}[htbp]
\begin{center}
\includegraphics[scale=0.2]{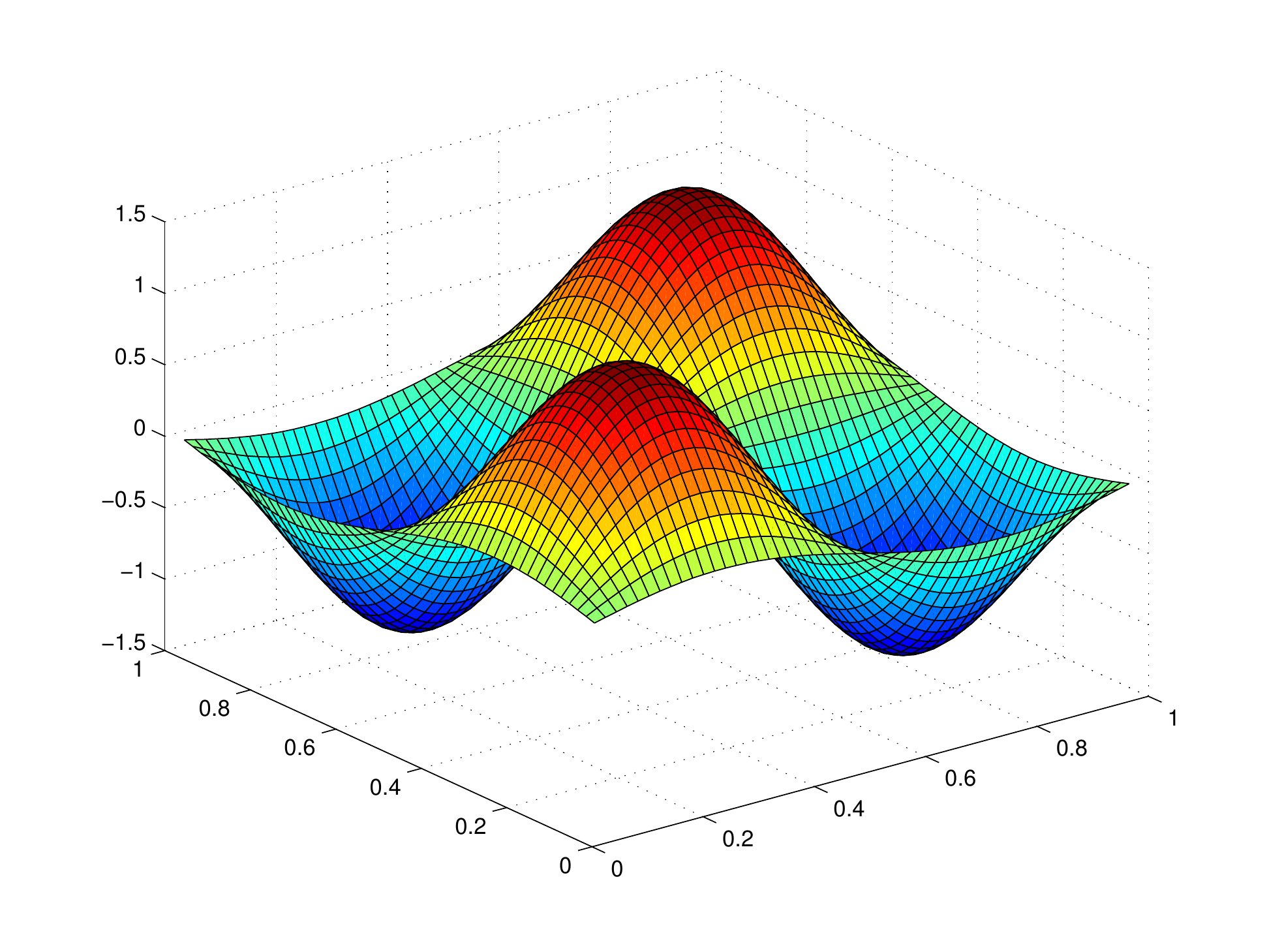}
\includegraphics[scale=0.2]{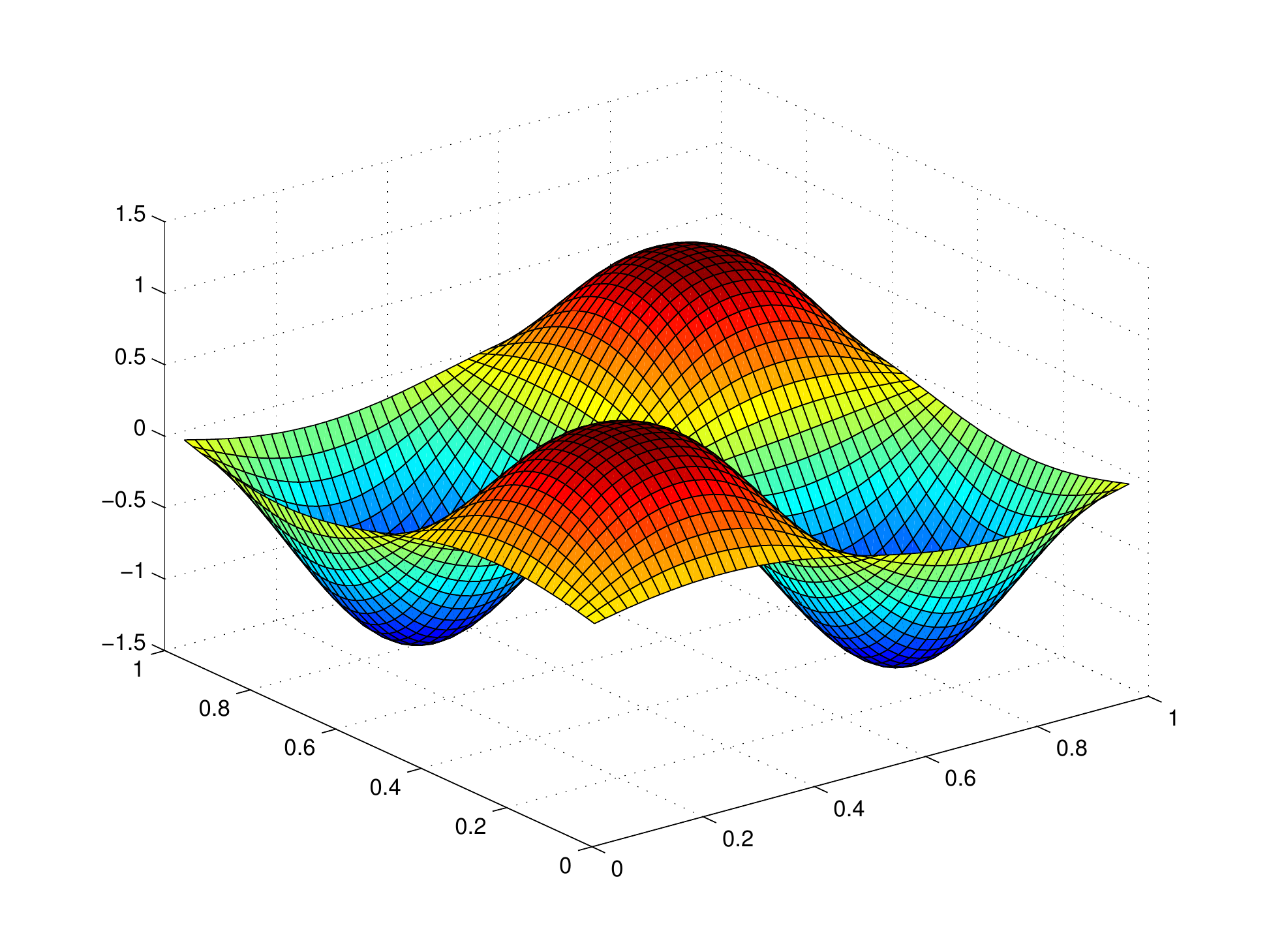}
\includegraphics[scale=0.2]{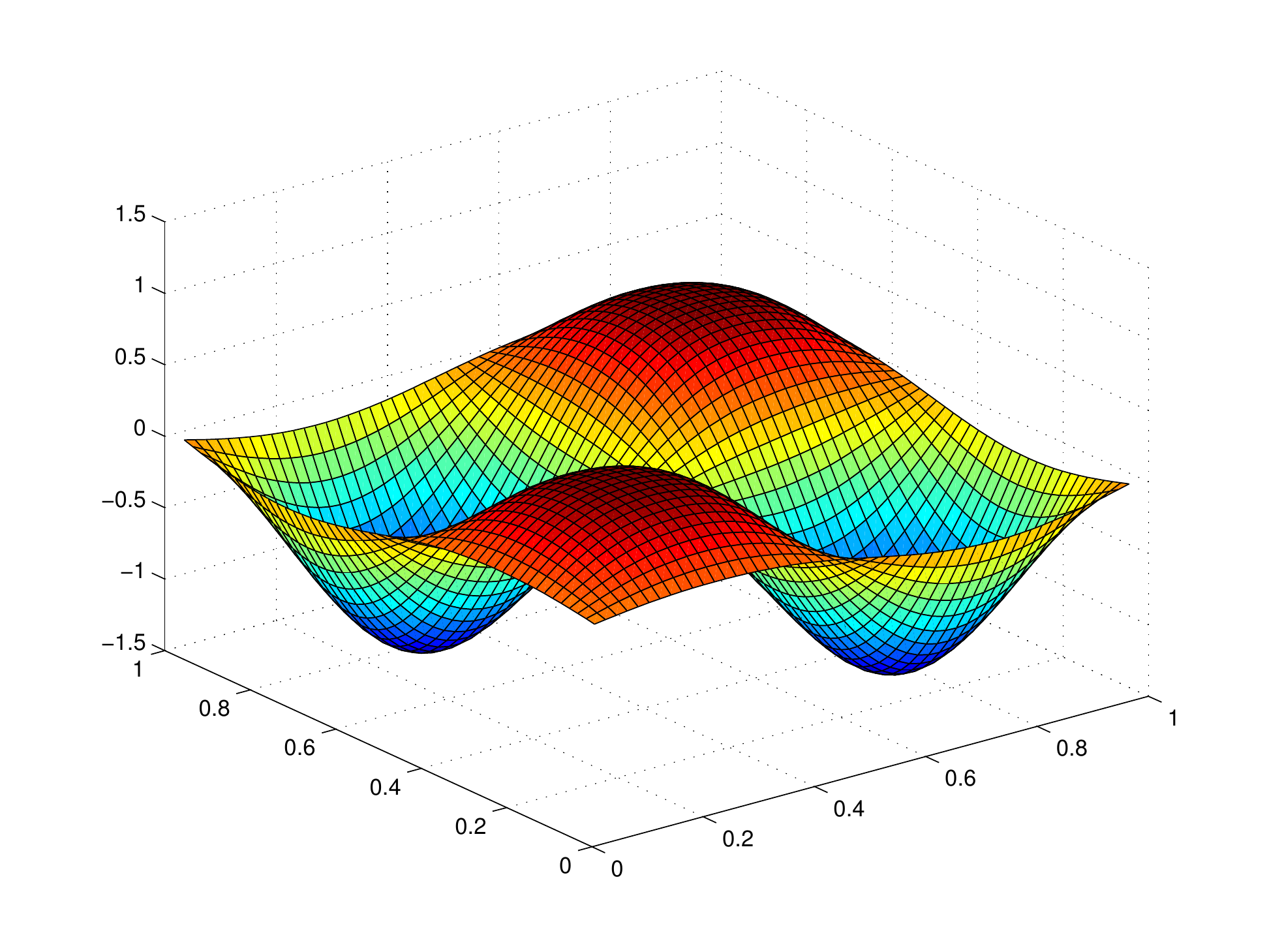}\\
\includegraphics[scale=0.2]{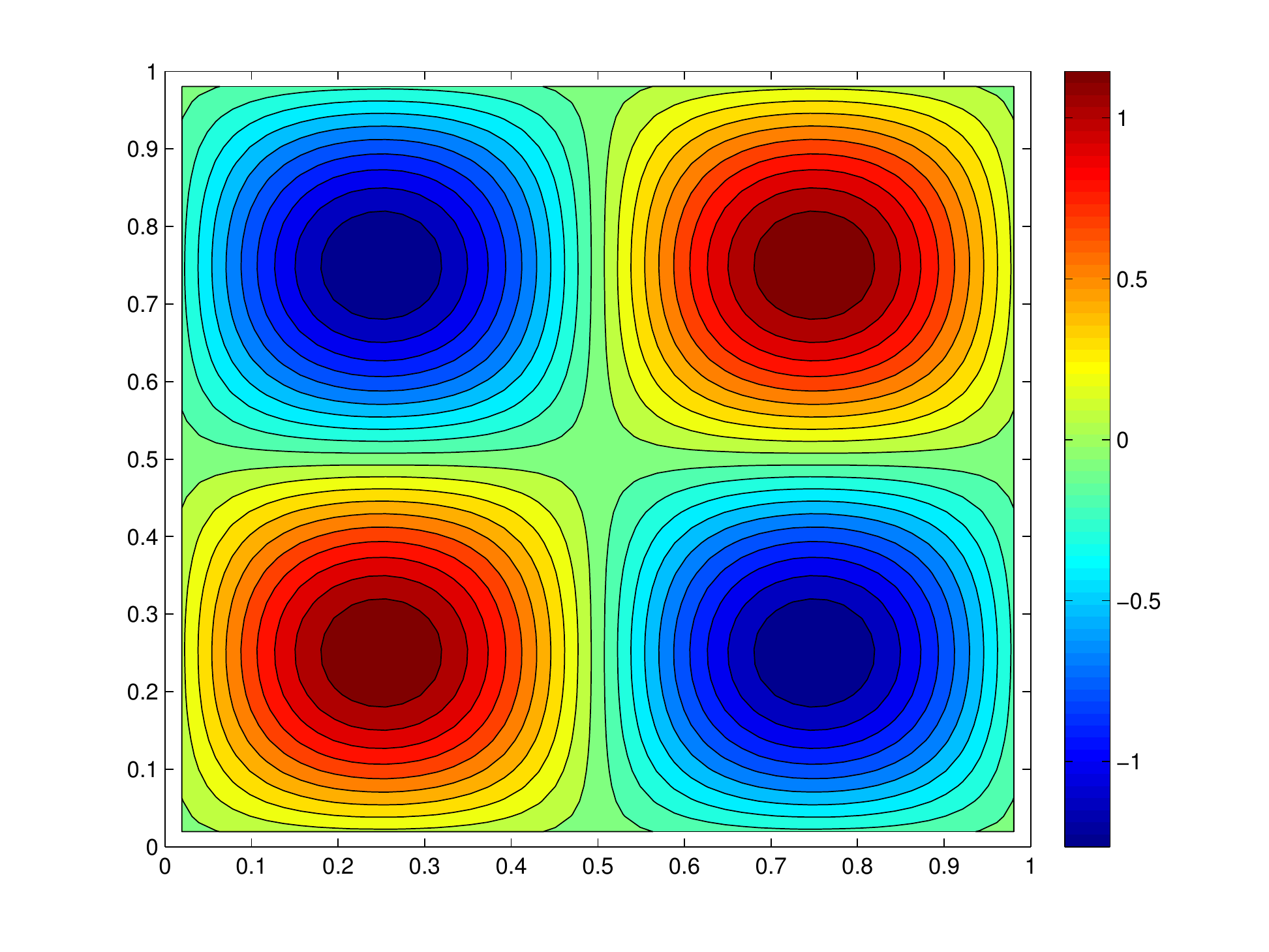}
\includegraphics[scale=0.2]{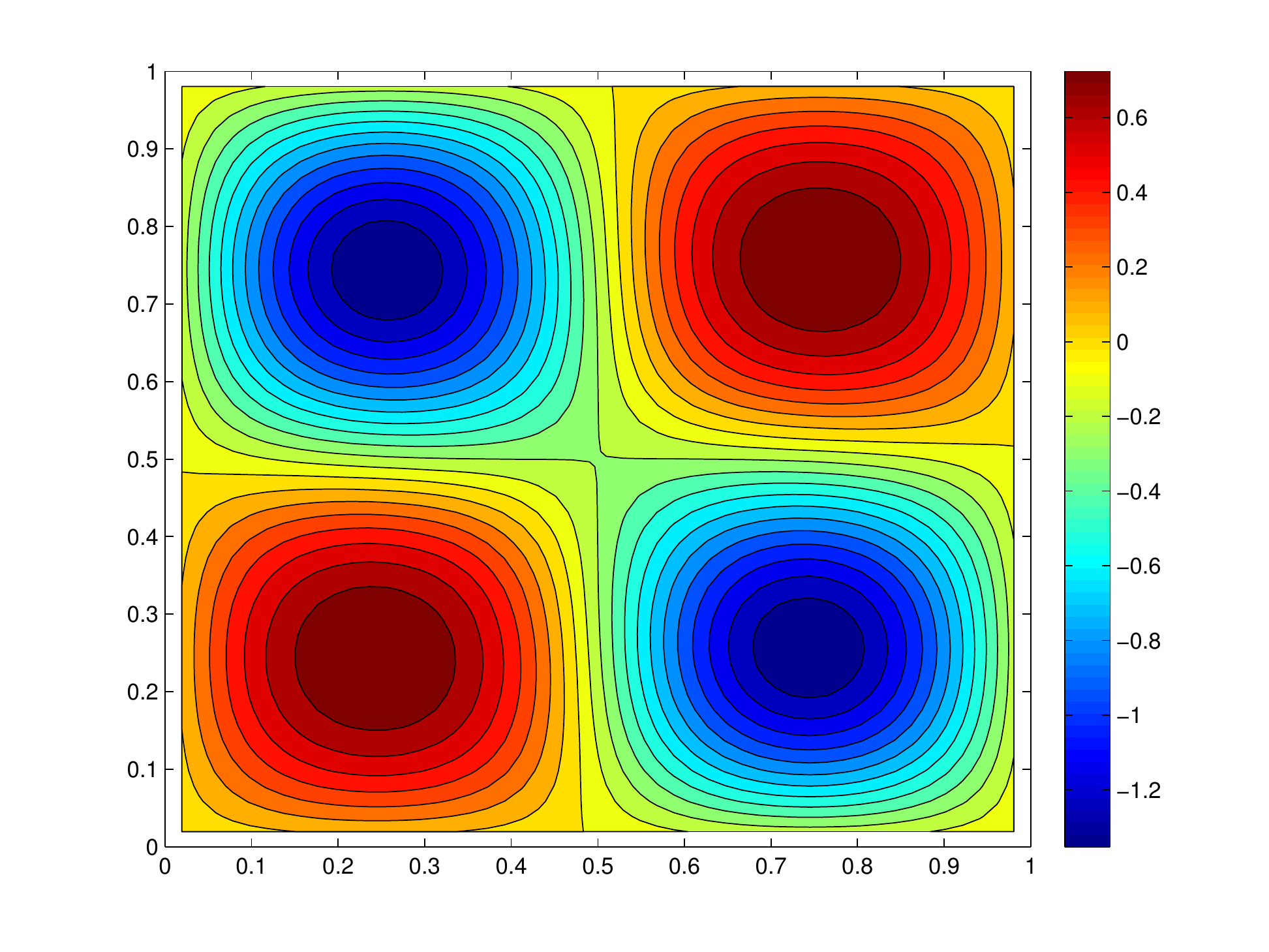}
\includegraphics[scale=0.2]{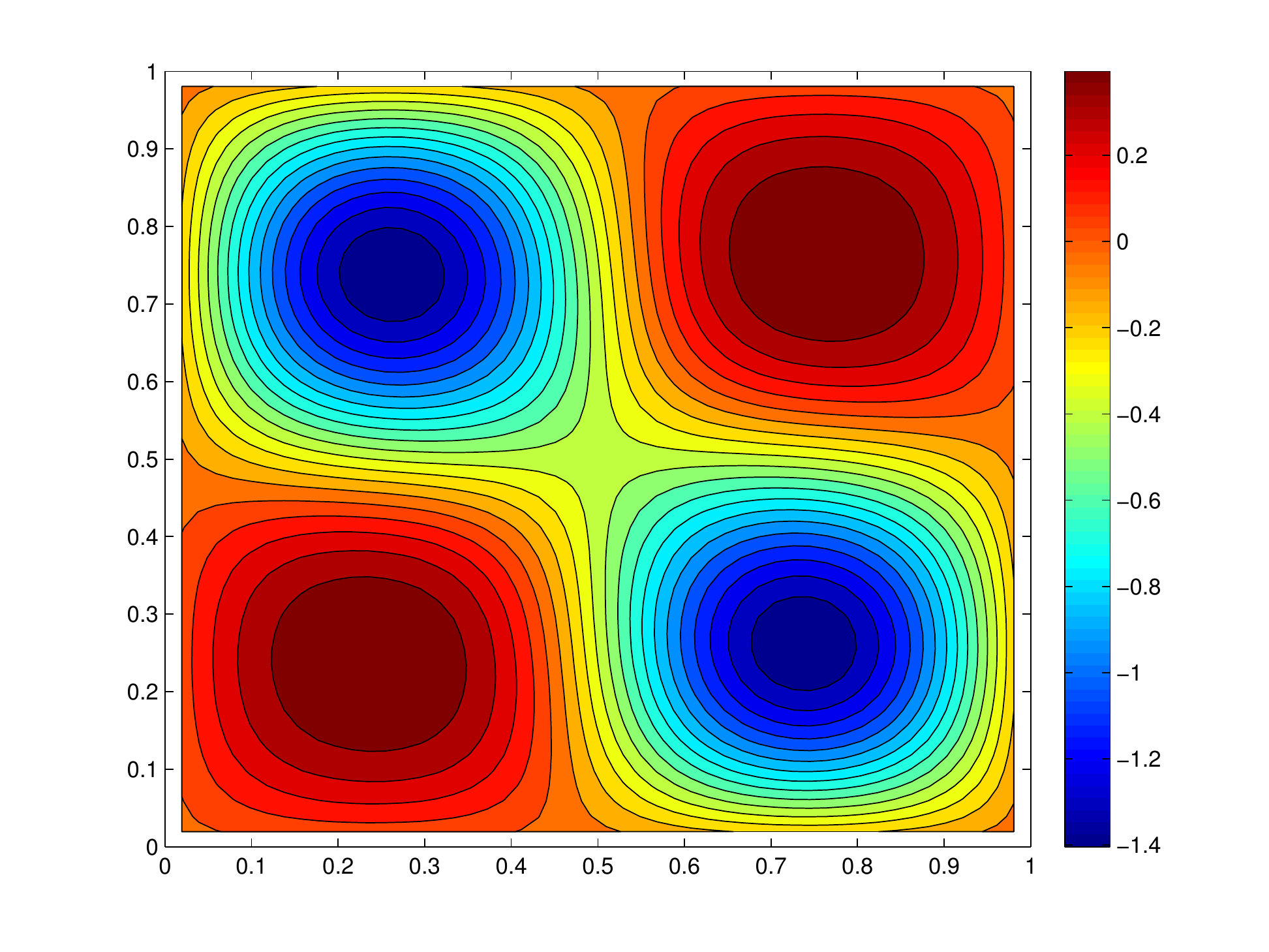}
\end{center}
\caption{Test 2: Solution of problem \eqref{test2:pde} for different parametric configurations (top) and count our lines (bottom)}
\label{test2:sol}
\end{figure}

The second numerical example concern a parametric elliptic equation. This example follows closely from \cite{CS10}. Let the dynamics given by:
\begin{equation}\label{test2:pde}
\left\{
\begin{aligned}
-\Delta u(x,y)+s(u(x,y);\mu)&=f(x,y) &\qquad (x,y)\in\Omega\\
u(x,y)&=0 &\qquad (x,y)\in\partial\Omega
\end{aligned}
\right.
\end{equation}
where the spatial variable $(x,y)\in\Omega=(0,1)^2$ and the parameters are $\mu=(\mu_1,\mu_2)\in\mathcal{D}=[0,01,10]^2\subset\mathbb{R}^2$ with a homogeneous Dirichlet boundary condition and nonlinearity
 $$s(u,\mu)=\dfrac{\mu_1}{\mu_2}(e^{\mu_2 u}-1),$$ and source term
 $$f(x,y)=100\sin(2\pi x)\sin(2\pi y).$$
 We numerically solve the system applying Newton's method to the nonlinear equations resulting from a FD discretization. The full dimension of the discretized problem is $n=2500$. 
% 
%Approximation with FD ($\Delta x=0.01=\Delta t, y_0(x)=x-x^2, \theta=1, \mu=11$)
The solution of \eqref{test2:pde} is shown in Figure \ref{test2:sol}.  Note that different choice of the parameter configuration leads different solutions.

\begin{figure}[htbp]
\begin{center}
\includegraphics[scale=0.25]{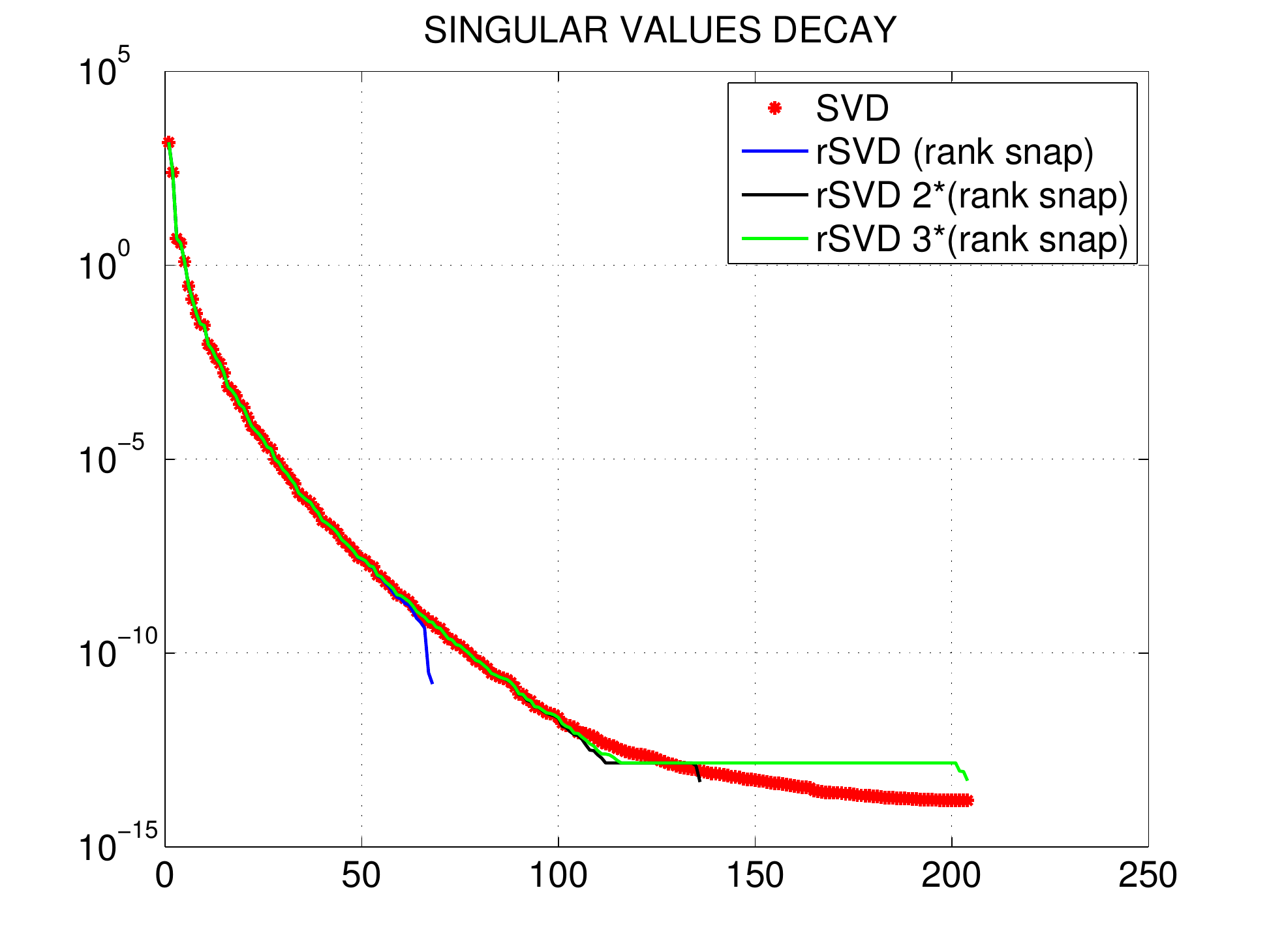}
\includegraphics[scale=0.25]{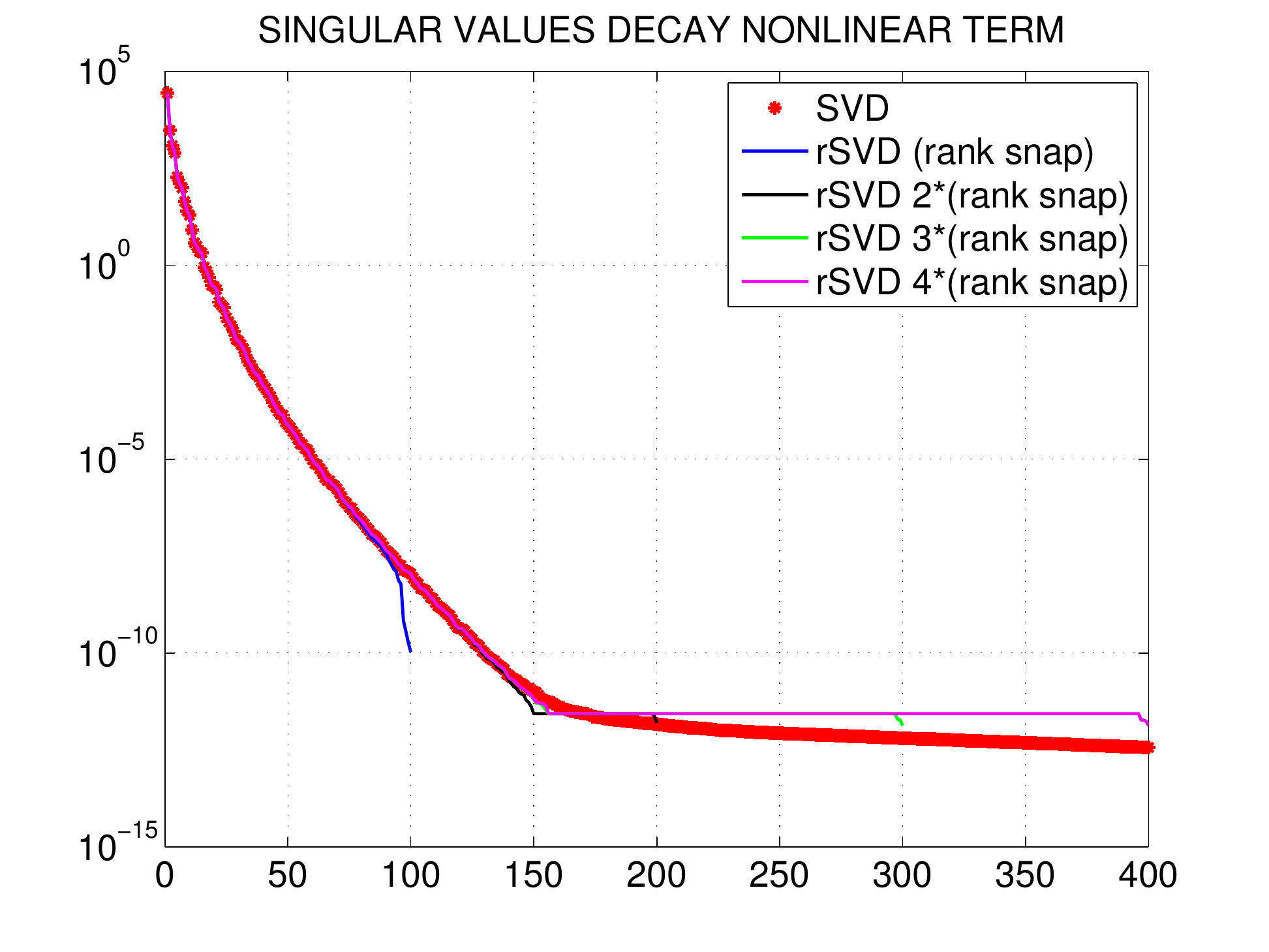}\\
\end{center}
\caption{Test 2: Decay of the singular values for the snapshot set with different number of measurements (left) and for the nonlinear term (right) }
\label{test2:svd}
\end{figure}

\begin{figure}[htbp]
\begin{center}
\includegraphics[scale=0.25]{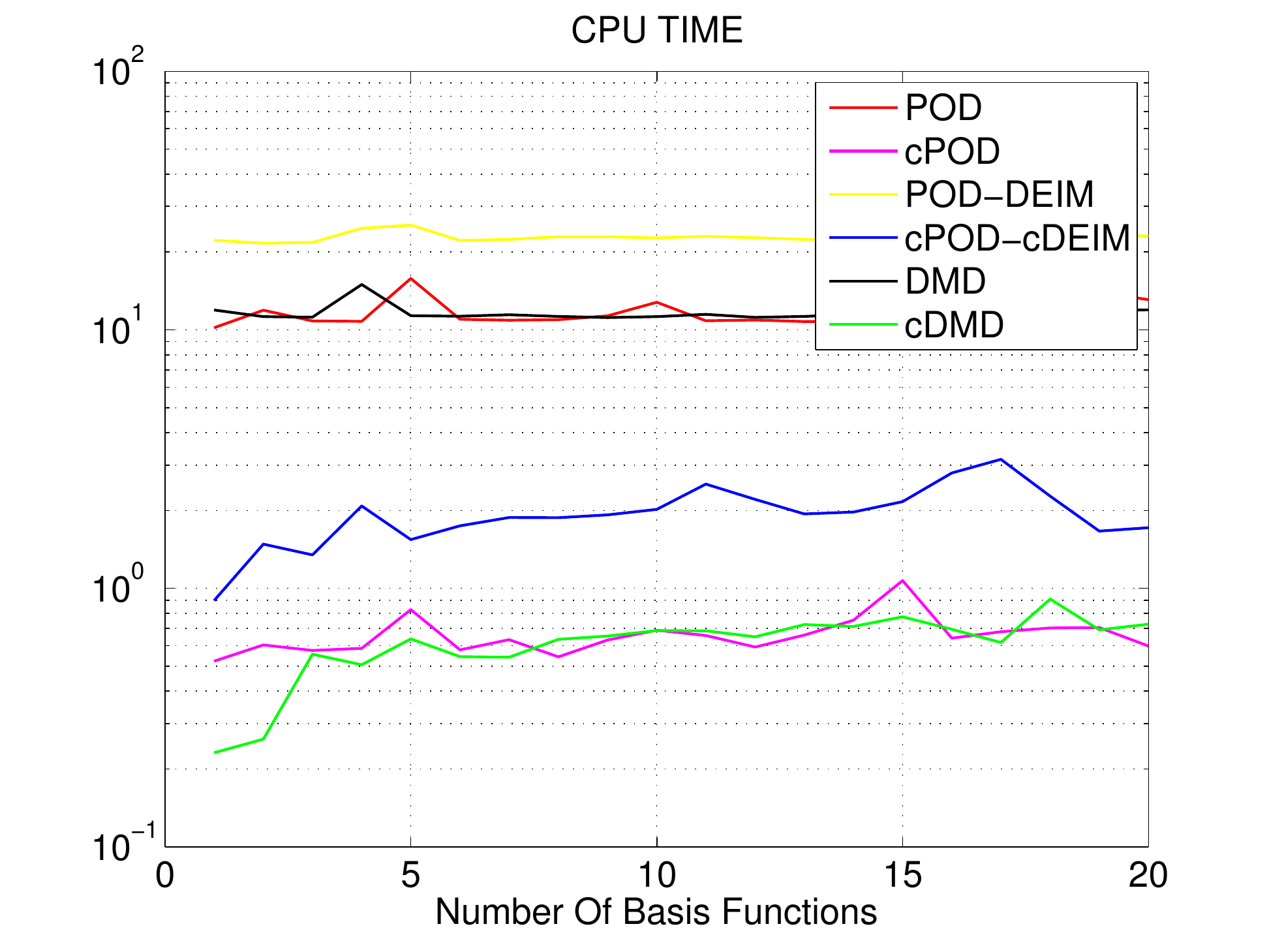}
\includegraphics[scale=0.25]{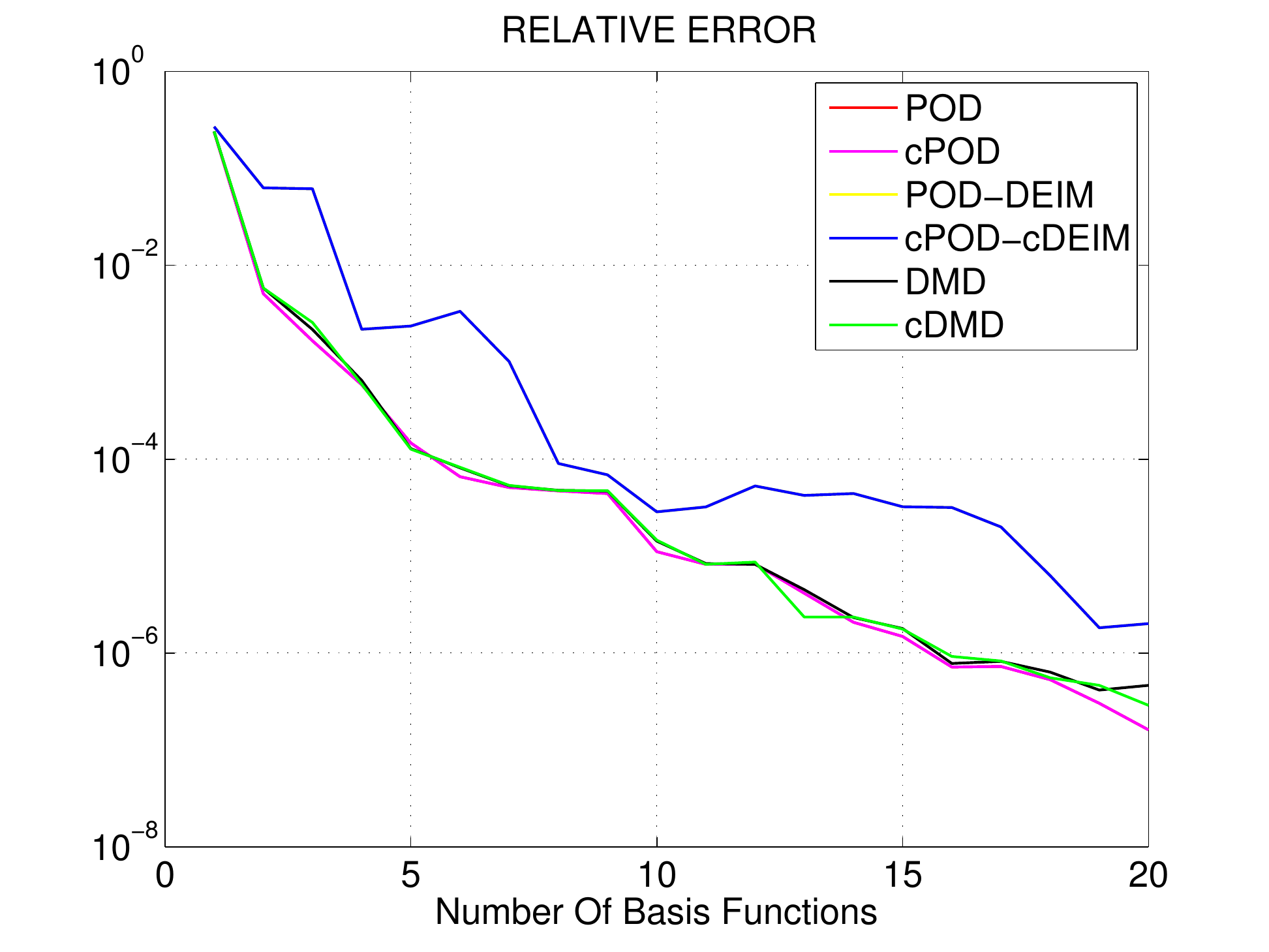}\\
\end{center}
\caption{Test 2: CPU-time of the offline-online stages (left) and Relative Error in Frobenius norm (right). We compare the following methods: POD (red), cPOD (magenta), POD-DEIM (yellow), cPOD-cDEIM (blue, DMD (black), cDMD (green). Number of model are always the same for all the methods.}

\label{test2:cpu}
\end{figure}

For the sake of completeness we show the decay of the singular values of the snapshot matrix on the left panel of Figure \ref{test2:svd}. We compare the singular values computed with the standard SVD and the randomized SVD as we increase the sampling points of the original matrix. As expected, it leads to improved approximations and, at the same time, faster approximations of the problem. A similar behavior comes from the nonlinear term (right panel).
Finally, we show the CPU time for all the methods studied in the left panel of Figure \ref{test2:cpu} and the error behavior in the right panel. They provide a similar analysis discussed in the previous test.

\section{Conclusion}
Model order reduction is a successful and commonly used technique that projects nonlinear high dimensional dynamical systems and PDEs into low dimensional surrogate models using optimal basis functions computed from information of the system.  Although the solution of the surrogate model is computationally efficient, the computation of the basis functions remains computationally expensive. In this paper we have demonstrated through several examples that compressed (randomized) techniques are a promising approach to circumventing expensive offline stages in model order reduction. In particular, when dealing with large snapshot matrices we suggest the use of randomized singular valued decomposition for the Proper Orthogonal Decomposition and compressed Dynamic Mode Decomposition. They both provide very accurate solutions and promise significant computational savings in the offline stage, which turns out to be the most expensive part of the building block for the surrogate model. 

Critical for enacting these computational enhancements is the advent of randomized linear algebra techniques.
Randomized linear algebra methods have been recently surveyed in \cite{gunnar_rLA}.  Indeed, the methods
are continuing to mature and have many critical error bounds associated with their proposed matrix factorizations.  
These efficient matrix decompositions are tremendously important for analyzing high dimensional data sets and/or for producing the low-dimensional subspaces required for ROMs.  Randomized techniques have continued to experience 
modifications that increase their efficiency and broaden the range of applicability of the methods. More broadly, 
randomized methods have application to classical (non-randomized) techniques for solving the same problems such as, e.g., Krylov methods, subspace iteration, and rank-revealing QR factorizations.   

Ultimately, ROMs are primarily concerned with producing rapid evaluation of surrogate models that represent the original high-dimensional system with a given accuracy.  Given the significant computational bottleneck for evaluating the low-dimensional projection, it is surprising the randomized linear algebra techniques have yet to penetrate the ROMs community.  We have explicitly demonstrated that such randomized techniques can be a significant enhancement of the ROMs architecture.  It should be used whenever possible given the current maturity of the technique and the error bounds available.

%\newpage

\end{document}